\documentclass[12pt, leqno]{article}
\usepackage{fullpage} 
\usepackage{amsmath}
\usepackage{amssymb}
\usepackage{natbib}
\usepackage{undertilde}
\usepackage{stmaryrd}
\usepackage{bbding}
\usepackage{enumerate}
\usepackage{graphicx}
\usepackage{xypic}
\xyoption{curve}
\usepackage{multicol}
\xyoption{color}
\usepackage{setspace}
\usepackage{boxedminipage}
\usepackage{nicefrac}
\usepackage{multicol}
\usepackage{amsthm}
  \usepackage[pdftex,
  bookmarks = false,
  pdfstartview = FitBH,
  linktocpage = true,
  pagebackref = true,
  pdfhighlight = /O,
  colorlinks=true,
  linkcolor=blue,
  citecolor=blue,
  filecolor = blue,
  urlcolor = blue,
  menucolor = blue,
]{hyperref}

\numberwithin{equation}{section}
\numberwithin{enumi}{section}
\numberwithin{foo}{section}
\numberwithin{thm}{section}

\numberwithin{equation}{section}
\numberwithin{enumi}{section}
\numberwithin{foo}{section}
\numberwithin{thm}{section}

\title{The Prehistory of the Subsystems of Second-Order Arithmetic}

 \author{Walter Dean\footnote{Department of Philosophy, University of Warwick, Coventry, CV4 7AL, United Kingdom, {\it E-mail}: W.H.Dean@warwick.ac.uk} \ and Sean Walsh\footnote{Department of Logic and Philosophy of Science, 5100 Social Science Plaza, University of California, Irvine, Irvine, CA 92697-5100, U.S.A., {\it E-mail}: swalsh108@gmail.com or walsh108@uci.edu}}

\begin{document}

\makeatletter
\renewcommand*{\@fnsymbol}[1]{\ifcase#1\or{*} \else{**}\fi}
\makeatother

\maketitle


\begin{abstract}
This paper presents a systematic study of the prehistory of the traditional subsystems of second-order arithmetic that feature prominently in the reverse mathematics program of Friedman and Simpson. We look in particular at: (i) the long arc from Poincar\'e to Feferman as concerns arithmetic definability and provability, (ii) the interplay between finitism and the formalization of analysis in the lecture notes and publications of Hilbert and Bernays, (iii) the uncertainty  as to the constructive status of principles equivalent to Weak K\"onig's Lemma, and (iv) the large-scale intellectual backdrop to arithmetical transfinite recursion in descriptive set theory and its effectivization by Borel, Lusin, Addison, and others.
\end{abstract}

\newpage

\tableofcontents

\newpage


\section{Introduction}\label{sec:intro}

The reverse mathematics project of Friedman and Simpson (\cite{friedman1975}, \cite{friedman1976}, \cite{simpson2009}) has been one of the more active areas in mathematical logic in recent decades. This program aims to calibrate the set existence principles implicit in ordinary mathematics by showing such principles to be equivalent to one of a small handful of traditional subsystems of second-order arithmetic. In increasing order of strength, the traditional systems are known as $\mathsf{RCA}_0$, $\mathsf{WKL}_0$, $\mathsf{ACA}_0$, $\mathsf{ATR}_0$, and $\mathsf{\Pi^1_1\mbox{-}CA}_0$.

The aim of this paper is to set out the history of these constituent systems of the reverse mathematics enterprise. Some parts of this history are well-known, and are part of the folklore which one learns when one learns this subject. This includes: Weyl's formulation of a system very much like $\mathsf{ACA}_0$ in his 1918 \emph{Das Kontinuum} (\cite{weyl1918a}), and the roots of $\mathsf{RCA}_0$ in the finitism of Hilbert and Bernay's \emph{Grundlagen} (\cite{hilbert1934}, \cite{hilbert1939}).

However, hitherto there has been no attempt to set out in a systematic manner how we get from Weyl's 1918 book to Friedman's 1974 address at the International Congress of Mathematicians (\cite{friedman1975}), or how ideas related to Hilbert and Bernay's finitism have come to find a second~life as a base system which so many logicians today regularly employ. As we will see, this is not merely a latter-day rekindling of earlier foundational projects, but rather is a continuous intellectual development which spans generations and interacts with nearly every area of mathematical logic.

We primarily envision this paper as a historical companion to the first chapter of Simpson's monograph \citeyearpar{simpson2009}. We follow Simpson by beginning with $\mathsf{ACA}_0$ in \S\ref{sec:ACA0}, and then moving to $\mathsf{RCA}_0$ in \S\ref{sec:RCA0}, then $\mathsf{WKL}_0$ in \S\ref{sec:wkl_0}, and reaching finally $\mathsf{ATR}_0$ in \S\ref{sec:ATR0}.\footnote{We treat $\mathsf{\Pi^1_1\mbox{-}CA}_0$ primarily in the footnotes, largely because it is a ``rounding out'' of the other systems.} This order of presentation also happily agrees in large measure with the historical order of the development of the systems of the reverse mathematics enterprise.

But this study should also appeal to historians and logicians who do not have a vested interest in reverse mathematics as such. This is for two related reasons. First, the history of the subsystems of second-order arithmetic is a veritable crossroads for various important ideas in the history of logic of the last century. For instance, the history that we set out below contains the history of the formalization of the comprehension schema for second-order logic, and the history of the idea of a basis~theorem in computability theory. Indeed, it would be hard to discuss the history of such ideas without thereby writing a history of the subsystems of second-order arithmetic.

Second, Friedman and others have often emphasized that mathematical logic can be seen as the study of a well-ordered chain of theories of increasing interpretability strength (cf. \cite{friedman2007a} \S{7}, \cite{feferman2000a} p. 428, \cite{koellner2009} p. 100), starting with weak theories of arithmetic and reaching up into the large cardinal hierarchy. In this chain, second-order arithmetic occupies a crucial transition~point between the number-theoretic and the set-theoretic. The origins of the subsystems thus ought be of considerable interest even to those whose own research falls in the great~space below or the great~space above second-order arithmetic.

One regrettable limitation of our study is that we focus exclusively on the best known of the traditional subsystems treated in \citep{simpson2009}.  We hence 
pay little attention to those systems which are based on choice principles, which are based on finitary forms of Ramsey's Theorem, or which are weaker than $\mathsf{RCA}_0$. We similarly are silent upon the development of constructive reverse mathematics (e.g. \cite{veldman2014}) and systems that use higher-types in their formalization (e.g. \cite{kohlenbach2008}).  Further, for reasons of space, we also adopt the pretense that there is no~gap in-between Friedman's 1974 talk (\cite{friedman1975}) and the appearance of the first edition of Simpson's monograph (\cite{simpson1999a}). Needless to say, much important work was done in these years which shape our contemporary understanding of reverse mathematics, to which our history does little justice.

We should finally mention that in electing to concentrate on  subsystems of second-order arithmetic themselves, our survey also fails to explicitly highlight what is often presented as the most characteristic feature of reverse mathematics -- i.e. the many instances of so-called \emph{reversals} by which theorems of classical mathematics have been proven equivalent to the various axiomatic principles which define the subsystems over the base theory $\mathsf{RCA}_0$. We will return to discuss such results briefly in the concluding section \S\ref{sec:conclusions} in regard to the received view that reversals serve to measure the set existence assumptions which are necessary and sufficient to prove classical theorems (e.g. \citep{simpson1988a}, \citep[I.1, I.9]{simpson2009}). 

By taking subsystems rather than reversals as our focus, our study is able to  provide a context for discussing the development in fields like computability theory and descriptive set theory which are more closely tied to the specific setting of second-order arithmetic. By the same token, were we to take reversals as our focus, then \emph{a priori} it seems that this would amount to a study of the general enterprise of proving equivalences between mathematical statements and set-theoretic statements over any base theory. This would be a broad and ultimately different topic,\footnote{For example, such a study should presumably also include equivalent formulations of the Axiom of Choice (e.g. \citep{rubin1963}), and Dedekind's set theoretic development of analysis leading to \citep{dedekind1888} (whose role as an anticipation of reverse mathematics has been highlighted by \citet{sieg2005}).}  and one which would be less directly tied to the origins of reverse mathematics as presented in Simpson's monograph \citep{simpson2009}.

It is perhaps inevitable that any exposition of the ``prehistory'' of an intellectual enterprise will take on a somewhat whiggish cast.  Nonetheless, we have attempted to approached the development of the particular axiomatic systems in question not from the perspective of rational reconstruction, but rather by attending to the details of the specific contexts which led to their isolation.  For instance, something that comes out in our history is the many starts and stops along the way-- such as the long pause between Weyl and Grzegorczyk discussed in \S\ref{subsec:GMKeffectiveanalysis}-- and the moments of hesitation -- such as over the constructive credentials of G\"odel's completeness theorem (to which Weak K\"onig's Lemma reverses) discussed in \S\ref{sec:basisthm}. Moreover, what ultimately emerges in the history of our topic is less the triumph of any one viewpoint than the development of a neutral framework which may serve to chart the relationships between many distinct positions in the foundations of mathematics.

Before we begin with the history of $\mathsf{ACA}_0$ in the next section, let us recall the basics of the axioms of the subsystems with which we shall be concerned. The axioms of full second-order Peano arithmetic or $\mathsf{Z}_2$ start with a finite set of axioms saying how addition and multiplication interact with zero and successor (the axioms of Robinson's $\mathsf{Q}$), and adds to them the induction axiom
\begin{equation*}
\forall \; X \; [X(0) \; \& \; \forall \; n \; X(n)\rightarrow X(S(n))]\rightarrow \forall \; n \; Xn
\end{equation*}
\noindent and the comprehension schema:
\begin{equation*}
\exists \; X \; \forall \; n \; (\varphi(n)\leftrightarrow Xn)
\end{equation*}
\noindent wherein $X$ does not appear free in $\varphi(n)$. As these formulas indicate, the distinction between numbers and sets of numbers is marked by using lower-case roman letters for the former and upper-case roman letters for the latter. The subsystem $\mathsf{ACA}_0$ is formed by restricting the comprehension schema to formulas $\varphi(n)$ that contain no second-order quantifiers. The system $\mathsf{RCA}_0$ is formed by (i) restricting comprehension further to formulas $\varphi(n)$ which are recursive, in that both it and its negation can be expressed by a $\Sigma^0_1$-formula, i.e. one starting with an existential quantifier over numbers and followed by only bounded quantifiers; and by (ii) replacing the induction axiom by the induction schema over $\Sigma^0_1$-formulas. The system $\mathsf{WKL}_0$ is formed from $\mathsf{RCA}_0$ by the addition of a statement to the effect that ``every infinite subtree of the full binary tree has an infinite path''; see \S\ref{sec:wkl_0} for more details. Finally, the system $\mathsf{ATR}_0$ is formed from $\mathsf{RCA}_0$ by the addition of a statement to the effect that ``one can do transfinite recursion along any well-order with an arithmetic operator''; see \S\ref{sec:ATR0} for more details. In all our notation and terminology, we are following the first chapter of Simpson's monograph \citeyearpar{simpson2009}.



\section{Arithmetical comprehension and related systems}\label{sec:ACA0}


\subsection{Russell, Poincar\'e, Zermelo, and Weyl on set existence}\label{subsec:RPZWsetexistence}

The history of $\mathsf{ACA}_0$ and related systems is closely tied to the history of predicativity. The meaning of this term was forged gradually in the writings of Russell and Poincar\'e, beginning in the years 1905-1906.  In 1905, Russell used the word \emph{predicative} to demarcate those formulas which determine sets from those which did not; and he took responding to the paradoxes to involve saying which formulas were predicative and which were not (\cite{russell1907a} p. 34). Poincar\'e immediately appropriated this vocabulary and in 1906 proposed identifying the \emph{im}predicative with that which contains a vicious circle, albeit without attempting to make precise the relevant sense of circularity (\cite{poincare1906a} pp. 307-308). Russell then concured and initially proposed understanding the type of vicious circularity at issue in terms of self-applicability (\cite{russell1906a} p. 634). This proposal eventually evolved into the \emph{vicious circle principle}, to which we shall return in \S\ref{sec:KWFpredprov}.

Hence in its first usage proposed by Poincar\'e and Russell in the years 1905-1906, \emph{predicative} carries two meanings: it serves as a sufficient and perhaps necessary condition for set-existence, and it is indicative of a certain lack of circularity in definitions. To the modern ear, these two things sound rather different in character, and so too did they to Zermelo, who in 1907-8 set out to simultaneously axiomatize set-theory and respond to the criticisms of his 1904 proof of the well-ordering theorem. Zermelo suggested that the appeal to the existence of infimums and supremums on the real line, such as in the standard proof of the fundamental theorem of algebra,  possessed the circular features at which Poincar\'e and Russell had gestured, so that an insistence on tying set-existence to a lack of circularity in definitions would make ``science impossible'' (\cite{zermelo1908a} p. 118). In 1909-1910, Poincar\'e conceded the force of this kind of objection but responded that in the specific case pointed out by Zermelo, one could proceed by taking infimums and supremums of rationals and thus transform the original proof into one which adheres to the predicative restrictions (cf. \cite{poincare1909b} p. 199, \cite{poincare1910a} p. 48).

In his own 1908 axiomatization, Zermelo formulated the axiom of separation by saying that any ``definite property'' of an antecedently given set determines a subset (\cite{zermelo1908b} pp. 263). In 1910, Weyl, like many at the time (cf. \cite{moore1982a} \S{3.3} pp. 160 ff, \cite{ebbinghaus2003a} \S{2} pp. 199 ff), expressed frustration with the lack of precision in this formulation, and proposed an alternative formulation, saying: ``by a definite relation shall be understood that which is defined through finitely many applications of suitably modified definition principles on the basis of the two relations of equality and membership'' (\cite{weyl1910a}, \cite{weyl1968a} Vol. 1 p. 304). It is not too hard to see in this the kernel of the contemporary inductive definition of well-formed formula.\footnote{See \cite{feferman1998} pp. 258--259 for a formalization of Weyl's 1910 notion of ``definite relation.''}

Thus in his 1918 book \emph{Das Kontinuum}, Weyl formulates a second-order system whose first-order objects are natural numbers. Saying that he was motivated by a desire to ``fix more precisely'' Zermelo's notion of definite property (\cite{weyl1918a} p. 36), Weyl again presents the inductive definition of a well-formed \emph{first-order formula} and then says that to each such specifiable ``property $E$ there corresponds a set $({\bf E})$''  (\cite{weyl1918a} pp. 13, 31-32). As Feferman put it: ``Weyl's main step, then, was to see what could be accomplished in analysis if one worked [\ldots] only with the principle of arithmetical definition'' (\cite{feferman1998} p. 54).

In \emph{Das Kontinuum}, Weyl then explicitly attends to the Zermelo-Poincar\'e question of whether classical results like the fundamental theorem of algebra hold in this more restricted setting. In effect, one of the capstones of Weyl's book is the demonstration that his formalization allows one to establish the fundamental theorem (\cite{weyl1918a} p. 64, \cite{feferman1964a} p. 7). But this comes only after having set up the rudiments of the real number system and continuous functions on it. Weyl formalizes real numbers as Dedekind cuts of rationals (\cite{weyl1918a} p. 51), and he proves the completeness of the real line in the sense that every Cauchy sequence converges (\cite{weyl1918a} pp. 57-58). However, he cautions that in general the existence of infimums and supremums of \emph{arbitrary} bounded subsets of real numbers is ``in no way secured'' in his model (\cite{weyl1918a} p. 60).

In his 1921 paper ``On the New Foundational Crisis in Mathematics,'' Weyl slightly strengthened this conclusion, saying that one has to ``abandon'' the least upper bound principle in general and that there is no way to save it (\cite{weyl1921a} pp. 47-48). Much the same was expressed in the third section of his 1926 essay (\cite{weyl1926a}), where it was placed immediately subsequent to a discussion of the antinomies. This made it seem to some, like \cite{holder1926a}, that Weyl was suggesting that classical analysis was somehow touched by paradox. In the last section of his essay, H\"older notes that if one is given a countable sequence of real numbers, there is a way of constructing its infimums and supremums along lines accepted by Weyl (cf. \cite{holder1926a} p. 246 ff). But this is precisely what Weyl showed in his 1918 book, and Weyl's reservations about infimums and supremums was simply that the kinds of restrictions he was advocating do not guarantee their existence for arbitrary bounded sets of reals, but only for predicatively definable bounded sets of reals.


\subsection{Grzegorczyk, Mostowski, and Kond\^o on effective analysis}\label{subsec:GMKeffectiveanalysis}

The topic of Weyl's 1918 book was only taken up again in the mid 1950s.\footnote{Such a gap in the development of the study of predicativity is highlighted by \citet{feferman2005a} pp. 601-603.  In two other instances, Feferman had indicated that Grzegorczyk and Kond\^o were the intellectual successors to Weyl's program (cf. \cite{feferman1964a} p. 8, \cite{feferman1998} p. 291). This subsection is an attempt to fill out more of the details of this succession.} The intervening years had seen the development of computability theory, descriptive set theory, and proof theory, and it was against the backdrop of this enriched landscape that people began to reevaluate the predicative perspective. Thus Grzegorczyk opens his 1955 essay by saying that ``The purpose of this paper is to give strict mathematical shape to some ideas expressed by H. Weyl in `Das Kontinuum'~'' (\cite{grzegorczyk1955a} p. 311). In the bulk of the paper, Grzegorczyk proceed by studying analysis in the minimal $\omega$-model of $\mathsf{ACA}_0$, so that all second-order quantifiers were restricted to the arithmetically definable subsets of natural numbers. However, late in the paper he introduced an axiomatic version, and after stating the full comprehension schema says that ``We admit in this scheme only those formulas $[\ldots]$ in which each quantifier bounds a variable of the lowest type'' (\cite{grzegorczyk1955a} pp. 337-338). 

Despite this brief aside on axiomatization, it seems safe to say that Grzegorczyk's primary concern was with $\omega$-models. In the 1957 Amsterdam conference on constructivity (cf. \cite{krajewski2007a} p. 6), Grzegorczyk wrote that what distinguishes the Polish approach to constructivity was that ``All methods of proof are allowed. The constructive tendency consists only in the narrowing of the field of mathematical entities considered in classical analysis'' (\cite{grzegorczyk1959a} p. 43). Mostowski's contribution to the volume is similar in spirit to Grzegorczyk's: he looks explicitly at the minimal $\omega$-model of $\mathsf{ACA}_0$ and notes that it ``coincides with the universe of the constructive analysis of Hermann Weyl'' (\cite{mostowski1959a} p. 183); and he too notes that there is a natural axiomatization close to $\mathsf{ACA}_0$ (\cite{mostowski1959a} p. 184). 

An important passage in Mostowski's article, which we find no analogue of in Grzegorczyk, is the following, in which Mostowski suggests the project of trying to find multiple $\omega$-models for the \emph{full} comprehension schema:
\begin{quote}
We consider therefore a class ${\bf K}$ of sets of integers and ask, whether the comprehension axiom $ (EX) (x) [x\in X \equiv G(x)]$ is satisfied in ${\bf K}$; this is equivalent to the problem, whether all the axioms of classical arithmetic are satisfied in ${\bf K}$. [\ldots] [\P] [\P] The general recursive, elementarily definable, and hyperarithmetic systems are thus essentially different from the classical arithmetic: if we take any of them as a basis for mathematics we must abandon some classical principles. This result which, of course, is not at all surprising suggests immediately the problem of finding as simple a model as possible for the classical axioms of arithmetic and of set theory. If such a model could be defined by constructivistic means we would have a constructivistic justification of the classical systems. (Models which we have in mind are absolute for integers, i.e. their integers are isomorphic with the ordinary ones). (\cite{mostowski1959a} pp. 186-187).
\end{quote}
\noindent In the last pages of the paper, he amplifies upon what he intends by ``constructivist justification.'' He says that the ``most promising feature [of constructivism] is that it wants to inquire into the nature of mathematical entities and to find a justification for the general laws which govern them, whereas platonism takes these laws as granted without any further discussion'' (\cite{mostowski1959a} p. 192). 

The sense of ``justification'' here is then more proximate to what we today might intend by ``explanation'': Mostowski seeks to identify some extension of the notion of computation (akin to the way \emph{hyperarithmetic} extends the notion of \emph{computable}) so that all of the comprehension axioms come out true when the second-order quantifiers are restricted to this extension. If, contrary to fact, all sets of natural numbers were computable in this more extended sense, one would then have an explanation for the truth of the full comprehension schema. As Motowski says, this is not constructivism in the sense of Brouwer, but ``gives merely a glance on constructivism, so to say, from the outside'' (\cite{mostowski1959a} p. 180). Heyting, in his contribution to the 1957 conference, concurred with this, saying that in a genuinely constructive theory ``there can be no mentioning of other than constructible objects.'' Hence, in Heyting's eyes, Mostowski and Grzegorczyk's preference for working in a classical metatheory precludes their work from being constructive in this more austere sense (\cite{heyting1959a} p. 69).

The work of this Polish school was heavily influenced by developments in computability theory, which during this period had tended to be closely related to developments in descriptive set theory. For instance, Kleene tells us that the identification of recursive sets of natural numbers with the $\Delta^0_1$-definable sets was first suggested by Souslin's theorem (cf. \cite{kleene1955a} p. 196). In the French tradition, work in descriptive set theory was associated to ideas of Poinca\'re and Lebesgue which sometimes go under the heading of French semi-intuitionism (cf. \cite{michel2008a}). This was the tradition in which Kond\^o worked, and on at least two occasions he begins his papers with invocations of the claims of Poincar\'e and Lebesgue that the only objects in mathematics are those which can be defined in a finite number of words (\cite{poincare1909a} p. 482, \cite{lebesgue1905a} p. 205, \cite{kondo1956b}, \cite{kondo1985a} p. 330, \cite{kondo1958a} p. 1). This requirement was practically implemented by restricting attention to either implicitly definable sets, or to classes of sets like the Borel or projective sets.

Kond\^o's approach in \citeyearpar{kondo1958a} was to consider two subfields $k_0\subseteq k$ of the real numbers, and to consider the model $\mathcal{A}(k_0,k)$ which consists of first-order part $k$ with a distinguished predicate for the integers and for $k_0$, and with second-order part consisting of all first-order definable subsets of $k$ in the ring signature expanded by these two predicates (cf. notation for polynomials~$F$ on \cite{kondo1958a} pp. 12-13, the operation $\mathcal{L}F$ on \cite{kondo1958a} pp. 13-14, and $\mathcal{A}(k_0,k)$ on \cite{kondo1958a} p. 19). The structures  $\mathcal{A}(k_0,k)$ are called, in Kond\^o's terminology, models of \emph{relative analysis}. In the case where $k_0$ is the rationals and $k$ is the reals, the second-order part of this model consists of the projective sets. In the case where both $k_0$ and $k$ are the rationals, this model will be a notational variant of the minimal $\omega$-model of $\mathsf{ACA}_0$. Kond\^o also considers ways to map one model of relative analysis to another. He considers the map $\mathcal{A}(k_0,k_0)\mapsto \mathcal{A}(k_0, \pi(k_0))$ where $\pi(k_0)$ consists of all the reals whose Dedekind cut is a set in $\mathcal{A}(k_0,k_0)$ (cf. \cite{kondo1958a} p. 74). He indicates in a later paper (\cite{kondo1960a} p. 62) that when $k_0$ is the rationals then the model of relative analysis $\mathcal{A}(k_0,\pi(k_0))$ is closely related to the systems of Weyl and Grzegorcyzk.\footnote{The mathematical result for which Kond\^o is now most well-known is the uniformization theorem for coanalytic sets (\cite{moschovakis2009a} p. 178, \cite{kechris1995a} p. 306). But he proved this result in 1939 (cf. \cite{kondo1939}), and it does not play a role in the papers from the late 1950s and early 1960s. Of course, retrospectively we can see a reason for this: the uniformization theorem reverses to $\mathsf{\Pi^1_1\mbox{-}CA}_0$ over $\mathsf{ATR}_0$ (cf. \cite{simpson2009} p. 225), and hence requires stronger assumptions than the perspective which Kond\^o was exploring in the papers from the late 1950s and early 1960s.}


\subsection{Kreisel on predicative definability}\label{subsec:Konpreddefn}

In his review of \cite{kondo1958a}, Kreisel noted that Kond\^o's various results about what is common to all models of relative analysis ``may also be expected on axiomatic grounds'' (\cite{kreisel1959a}). This idea of connecting \emph{multiple models} to an \emph{axiomatic treatment} is reiterated in Kreisel's retrospective 1976 article \citep{kreisel1976a}, where he notes that ``the bulk of current theorems \emph{generalize}: wherever the notion of \emph{set} is used, explicitly or implicitly, it may be interpreted to mean: set of a (so to speak `elementary') collection of sets satisfying the particular `weak' closure conditions'' (\cite{kreisel1976a} p. 109). 

Part of the impetus for this work came out of Kreisel's writings on the Hilbert Program. Kreisel noted in 1958 that the effects of the G\"odel incompleteness theorem for the Hilbert Program only touch formalized notions of provability, which might be different than the ``absolute'' notions of provability associated with traditional programs like finitism, constructivity, and predicativism (\cite{kreisel1958a} p. 177, cf. \cite{kreisel1968b} p. 323). However, Kreisel noted in a paper appearing two years later that ``predicative provability'' might be a rather different thing than ``predicative definability'' (\cite{kreisel1960a} p. 298, cf. \cite{wang1974a} p. 128). As an initial suggestion for how to formalize the latter, he proposed that one call a subsystem of second-order arithmetic \emph{predicative} if it has a unique minimal $\omega$-model (\cite{kreisel1960a} p. 298).

The idea then was to allow for a broader notion of predicativity than that which appears in Weyl or Grzegorcyzk. For, at around the same time as the first paper of Grzegorcyzk, it had been shown by Kleene that the $\Delta^1_1$-definable sets of natural numbers formed the minimal $\omega$-model of $\Delta^1_1$-comprehension (\cite{kleene1955a}). Kreisel wrote of Kleene's work that it ``provides a precise and satisfactory definition of the notion of predicative sets (of integers)'' (\cite{kreisel1955c}).  Further, Kreisel, Gandy and Tait subsequently proved that the $\Delta^1_1$-definable sets are the \emph{only} sets that are contained in the intersection of all $\omega$-models of a given recursively enumerable $\omega$-consistent set of axioms (in the language of second-order arithmetic).\footnote{Indeed, their result is much stronger in that one can replace `recursively enumerable' by `$\Pi^1_1$' (\cite{gandy1960} p. 579, cf. \cite{apt1973} p. 188).}

These notions deserved to be called \emph{predicative}, in Kreisel's view, because when one looks back at the original works of Poincar\'e and Weyl, there is a stability idea. For instance, Poincar\'e tells us that a predicative classification is one which is not changed by the introduction of new elements (\cite{poincare1909a} p. 463, cf. \cite{walsh2015a} \S{4}). This could be formalized, Kreisel noted, with the contemporary model-theoretic notion of absoluteness: a formula $\varphi(\overline{x})$ is \emph{absolute} between a substructure $\mathcal{M}$ and a superstructure $\mathcal{N}$ if $\mathcal{M}\models \varphi(\overline{a})$ iff $\mathcal{N}\models \varphi(\overline{a})$ for all tuples $\overline{a}$ from the substructure $\mathcal{M}$. Of course, formulas which are provably $\Delta^1_1$ in a subsystem of arithmetic are absolute between models of that subsystem which share the same first-order part (\cite{kreisel1960d} p. 378, \cite{kreisel1970c} p. 512, cf. \cite{feferman1987} p. 450). 

Connecting this back to the idea of a common set of axioms with multiple models, Kreisel wrote that the idea was to
\begin{quote}
find \emph{a convenient set of axioms for second or higher order arithmetic which are valid both when the variables of higher type} (in the sense of the \emph{simple} theory of types) \emph{are interpreted as ranging over all sets} (of the type considered) \emph{and when they are interpreted as ranging over predicative sets} (\cite{kreisel1962b} p. 311).
\end{quote}
\noindent It was in this paper that Kreisel formulated the $\Sigma^1_1$-choice axiom which forms the backbone of the subsystem which we now call $\mathsf{\Sigma^1_1}\mbox{-}\mathsf{AC}_0$ (\cite{kreisel1962b} p. 313).


\subsection{Kreisel, Wang, and Feferman on predicative provability}\label{sec:KWFpredprov}

Another related conception of predicativity pertained to predicative conceptions of proof. This was mentioned briefly in the previous section, but its roots go back to Russell's type theory. As is well-known, Russell did not begin with just simple type theory but rather with so-called \emph{ramified type theory}, which introduced different layers of higher-order variables. This first appeared in Russell's papers \citep{russell1908a,russell1910a} and was employed in the \emph{Principia} \citep{whitehead1910a}. In such systems, not only could one define a first round of second-order objects by first-order comprehension, but one could then define a second round of second-order objects by quantifying over first-order objects or second-order objects of the first round; and there are similarly $\omega$-many rounds for the third-order objects, the fourth-order objects, etc. (cf. \cite{schutte1960a} \S{27} pp. 245 ff, \cite{schutte1977a} \S{22} pp. 197 ff, \cite{church1976a}, \cite{hazen1983a} pp. 343~ff, \cite{urquhart2003a} \S{4} pp. 293 ff). If, as is common, one omits all but the first- and second-order objects, then the system is usually called \emph{ramified analysis}.

Both the mathematical elegance and philosophical motivation of ramified type-theory and ramified analysis have been subject to dispute. They were thought to be inelegant because the typing of second-order variables had no analogue in the mathematics which one sought to formalize in the system. In the words of Feferman, it was a ``parody of classical analysis'' (\cite{feferman1964a} p. 12, cf. \cite{kreisel1962b} p. 68, \cite{ramsey1926a} p. 186). But the rejoinder always was, as Wang once put it, that ``in formalizing actual proofs we do not have to let even the distinction of orders intrude'' (\cite{wang1954a} p. 266). 

  As for its philosophical motivation, out of his debate with Poincar\'e (cf. \cite{russell1906a} p. 634), Russell eventually settled upon the following formulation of the vicious circle principle: ``If, provided a certain collection had a total, it would have members only definable in terms of that total, then the said collection has no total'' (\cite{whitehead1910a} p. 40, \cite{whitehead1962a} p. 38). Ramsey noted that common everyday uses of definite descriptions like ``the tallest man in the room'' seem to violate the this principle (\cite{ramsey1925a} p. 368). Without trying to meet this objection directly, G\"odel in 1944 suggested that the principle might be seen to at least \emph{follow} from a ``constructivist'' conception of properties, on which they are built up iteratively out of definable sets (cf. \citep{godel1990} p. 127 and the footnote on  p. 119).\footnote{See \cite{parsons2002} \S{5} for a more thorough discussion of G\"odel's views. More generally, \cite{parsons2002} is a study of figures such as Hilbert, Bernays, and Ramsey as critics of ``definitionalism'', the view that ``sets are definable sets'' (\cite{parsons2002} p. 386).}

In his 1954 paper, \cite{wang1954a} described a ramified type theory which contained levels corresponding to infinite ordinals. However, Wang himself admitted that he did not really know how far up the ordinal hierarchy the system went (\cite{wang1954a} pp. 247-248, pp. 260-261, \cite{wang1955a} \S{9} pp. 77 ff). Kreisel then suggested only including the \emph{predicative ordinals}. Kreisel's proposed definition of predicative ordinal was inductive, and read as follows: if $\alpha$ is a predicative ordinal and ramified analysis up to level $\alpha$ proves that a recursive well-order, with order type $\beta$, is a well-order then $\beta$ is a predicative ordinal (cf. \cite{kreisel1960a} \S{5} p. 297, \cite{feferman2005a} p. 607, \cite{pohlers1987a} p. 413). Feferman and Sch\"utte  independently showed in \cite{feferman1964a}, \cite{schutte1965b, schutte1965a} that the least non-predicative ordinal is the ordinal which now bears the name of the \emph{Feferman-Sh\"utte} ordinal (cf. \cite{feferman1998} p. 121--122).

\subsection{$\mathsf{ACA}_0$ as a formal system}\label{sec:aca0formal}

We have just seen that some of the the first systems to attempt to axiomatize predicative reasoning included either ramified comprehension or iterations of provability along well-orders. Neither of these approaches leads directly to a unique characterization of the system $\mathsf{ACA}_0$.  

One development which anticipates the isolation of this system more directly was work on what is now known as G\"odel-Bernays set theory $\mathsf{GB}$.  This system was originally proposed as a two-sorted first-order theory of sets and classes by \citet{bernays1937}, based on a prior axiomatization by \citet{von-neumann1925}.  In addition to various axioms of set existence -- e.g. Empty Set, Pairing -- $\mathsf{GB}$ contains axioms formalizing the so-called G\"odel operations. These assert that the domain of classes is closed under operations such as complementation, intersection, and taking converses (in the case of classes which are relations).  Although this system is finitely axiomatized, \citet[p. 72]{bernays1937} showed that these operations are sufficient to prove the existence of all classes which are definable by formulas not containing bound class variables. On this basis \citet{mostowski1950} showed that $\mathsf{GB}$ is a conservative extension of $\mathsf{ZF}$.  

It is now known that $\mathsf{ACA}_0$ shares many of these features with $\mathsf{GB}$.  For instance $\mathsf{ACA}_0$ is conservative over $\mathsf{PA}$. And although this system is typically presented as consisting of $\mathsf{Q}$ together with the arithmetical comprehension scheme and the induction axiom, $\mathsf{ACA}_0$ may also be finitely axiomatized on the basis of appropriately chosen variants of the G\"odel operations (see \citep[p. 154]{hajek1998}).  \cite{bernays1942} also observed that it was possible to formalize portions of analysis (inclusive of the existence of least upper bounds) in the system $\mathsf{S}$ consisting of $\mathsf{GB}$ without the Axiom of Infinity for sets, in a manner which is similar to the formalization in second-order arithmetic carried out by \citet{hilbert1939} (as discussed below).  While such a development relies on the standard formalization of arithmetic in set theory, it is also possible to interpret $\mathsf{ZF}$ together with the negation of the Axiom of Infinity in $\mathsf{PA}$ via the Ackermann interpretation (cf. \cite{ackermann1937}) -- a fact which was systematically exploited by Wang in his investigation of $\mathsf{S}$ and similar fragments of $\mathsf{GB}$ which he referred as ``predicative set theory'' (cf., e.g., \cite{wang1953}).  \citet[p. 184]{mostowski1959a} similarly observed that the ``the part of Bernays' axiomatic system of set theory which deals with construction of classes represents an axiomatization of a constructivistic (elementary [i.e. arithmetically] definable) notion of a set.''\footnote{Since the minimal $\omega$-model of $\mathsf{S}$ corresponds to the hereditary finite sets together with the predicatively definable classes thereof,  such an observation can (as Mostowski observed) be understood semantically in terms of the relationship between this structure and the minimal $\omega$-model of $\mathsf{ACA}_0$.  And although this appears not have been noted at the time, it can be understood proof theoretically since the mutual interpretability of $\mathsf{PA}$ and $\mathsf{ZF} - \mathsf{Infinity} + \neg \mathsf{Infinity}$ extends to show that of  $\mathsf{ACA}_0$ and the $\mathsf{S} + \neg \mathsf{Infinity}$.}   

A final antecedent for the isolation of $\mathsf{ACA}_0$ is provided by work in computability theory which was inspired by the arithmetization of G\"odel's completeness and incompleteness theorems -- results which we will see below also played an important role in the delineation of $\mathsf{WKL}_0$.  For on the one hand, $\mathsf{WKL}_0$ can be characterized as an extension of $\mathsf{RCA}_0$ in virtue of the fact that any $\omega$-model of this theory itself contains a countable coded $\omega$-model of $\mathsf{WKL}_0$ (see \citet[VIII.2.7]{simpson2009}).  But on the other hand, \citet{simpson1973} suggests that $\mathsf{ACA}_0$ can also be characterized as the weakest such extension for which it may be shown that if a recursively axiomatizable theory $\mathsf{T}$ itself possesses an $\omega$-model, then there is also an $\omega$-model of $\mathsf{T}$ plus the formalization of the statement ``there does not exist a coded $\omega$-model of $\mathsf{T}$''.\footnote{An antecedent to this result was first obtained as a corollary to G\"odel's second incompleteness theorem by \citet{rosser1950} (see also \citep{mostowski1956} and \citep[VIII.5.6]{simpson2009}). In the form just stated, however, the result was obtained by \cite{steel1975} as a corollary of the following computability theoretic fact: there is no arithmetically definable relation $P \subseteq 2^{\omega} \times 2^{\omega}$ which defines an infinite descending sequence of Turing degrees -- i.e. $\langle A_i\ : \ i \in \omega \rangle$ and $A'_{i+1} \leq_T A_i$ for all $i$, where $A_{i+1}$ is the unique set such that $P(A_i,A_{i+1})$.}


\section{Hilbert and Bernays, the \textsl{Grundlagen der Mathematik}, and recursive comprehension}\label{sec:RCA0}

Two historical claims frequently made in regard to reverse mathematics are that the study of  second-order arithmetic and it subsystems can be traced to Hilbert and Bernays's \textsl{Grundlagen der Mathematik} \citeyearpar{hilbert1934}, \citeyearpar{hilbert1939} and that the system now known as $\mathsf{RCA}_0$ is somehow related to what they describe as the \textsl{finiten Standpunkt}  and thus also more generally to the view known as \textit{finitism} (e.g. \cite{hilbert1922}, \cite{hilbert1926}) (e.g. \citep{simpson1988a}, \citep[\S I, IX.3]{simpson2009}, \cite{feferman1993}). The first of these claims pertains to the general logical framework employed in the \textsl{Grundlagen} while the latter pertains to a specific set of arithmetical axioms.  And although a close reading of this work lends some support to both contentions, one of our aims in this section will be to bring out some complexities in the conventional narrative relating both to the development of second-order logic by Hilbert and others during the 1910s-1940s and the use of formal systems like $\mathsf{RCA}_0$ to provide a precise delineation of the \textsl{finiten Standpunkt}.  


\subsection{From the Axiom of Reducibility to second-order arithmetic}\label{sec:fromAR2SOL}

Both Hilbert and Bernays's lecture notes from 1917-1923 \citep{hilbert2013a}, as well as Hilbert and Ackermann's textbook \textsl{Grundz\"uge der theoretischen Logik} \citep{hilbert1928} culminate in a discussion of Russell and Whitehead's system from the \emph{Principia Mathematica} \citep{whitehead1910a}.\footnote{See also \citep{mancosu2003} for more on the reception of \textit{Principia} by Hilbert and his collaborators.} Prior to this discussion, each of these texts dramatically advance upon the \emph{Principia} in that they isolate and study the fragment of this system corresponding to propositional logic and first-order predicate logic. Hence, unlike the \emph{Principia} itself, both these lecture notes and the Hilbert-Ackermann monograph are immediately accessible to the modern reader as they follow our contemporary way of introducing logic.

For our purposes, the crucial idea in both sources is the measured dissatisfaction with the Axiom of Reducibility and the ramified type-theory of the \emph{Principia}. As mentioned above in \S\ref{sec:KWFpredprov}, the idea of ramification was to define a first round of second-order objects by first-order comprehension, and then to define a second round of second-order objects by quantifying over first-order objects and second-order objects of the first round, and then continuing onto further rounds. Russell and Whitehead further articulated the Axiom of Reducibility, which postulated that everything obtained at the second round (or a later round) was in fact extensional with something obtained at the first round (\cite{whitehead1910a} vol. 1 pp. 58~ff, 161~ff, \cite{whitehead1962a} pp. 55ff, pp. 166 ff).

The primary concern with the Axiom of Reducibility, evinced by Hilbert and his collaborators, was that it vitiated the intended interpretation of the ramified system. For they conceived this interpretation to be one on which one started with a given structure and added on its collection of first-order definable subsets, and then added a second collection of subsets definable in a first-order way from those, etc. That is, the intended interpretation is close to the``constructivist'' conception mentioned by G\"odel in 1944 (cf. \S\ref{sec:KWFpredprov}).

In the lecture notes, Hilbert and his collaborators go onto note that: ``in an arbitrary choice of the primitive properties and relations [of the structure] one cannot in general be sure that the Axiom of Reducibility is satisfied'' (\cite{hilbert2013a} p. 487). Rather, the interpretation on which this axiom is definitely satisfied is one on which the entities in the first round are ``considered as something existing in and of themselves, so that its plurality does not depend on actually given definitions nor at all on the possibility of us giving a definition'' (\cite{hilbert2013a} p. 206, cf. p. 487, cf. \cite{parsons2002} p. 378).\footnote{Using contemporary terminology, one might say that the models in which the Axiom of Reducibility is definitely satisfied are those where the second-order entities of the first round, the second round, the third round, etc. are all provided by the application of the powerset operation to the underlying first-order domain.}

This situation, they suggest, leads to the following dilemma regarding ramified type theory:
\begin{quote}
[\ldots] either [(a)] it is handled purely formally, in which case it is imprecise and offers no guarantee of being without contradictions, or [(b)] the logical operations will be made precise contentfully (\emph{inhaltlich}) so that contradictions are excluded, but that one does not obtain the usual proof methods of analysis and set theory
(\cite{hilbert2013a} p. 488).
\end{quote}
\noindent To reverse the order of the dilemma, the thought seems to be that either (b) ramified type theory is taken without the Axiom of Reducibility and is thus \emph{inhaltlich} but does not succeed in obtaining analysis, or (a) ramified type theory is taken with the Axiom of Reducibility and is thus not \emph{inhaltlich} and offers no guarantee of consistency and needs be treated purely formally. The reason they think that the Axiom of Reducibility is necessary for analysis, in the context of ramified type theory, is that real numbers were being formalized as left Dedekind cuts, so that the least upper bound of a formula defining a bounded set of reals would be given with a second-order existential quantifier corresponding to the union of all the cuts. And this higher-order quantifier ought to range over all second-order objects and not just those from the first round (cf. \cite{hilbert2013a} pp. 213, 485, 906, \cite{hilbert1928} p. 111, \cite{hilbert1939} p. 463).

In subsequent writings, Hilbert and his collaborators clearly opted for horn~(a) of the dilemma. Since they preferred a system which was being treated purely formally and judged by the extent to which it was able to recover analysis, it was noted in the Hilbert-Ackermann monograph that one could simply remove the ramified system entirely and move to what we would now call simple type theory (\cite{hilbert1928} p. 115, \cite{hilbert2013a} p. 909). And in the second edition of the Hilbert-Ackermann monograph in 1938 one finds the following statement of the comprehension schema:
\begin{quote}
{}Let $G_1, G_2, \ldots, G_n$ variables of any type $a_1, \ldots, a_n$, and $F$ a variable of type $(a_1, \ldots, a_n)$, and $A(G_1, \ldots, G_n)$ a formula that has free variables $G_1, G_2, \ldots, G_n$. Then each formula of the form [\P] (V) $(EF) (G_1) \ldots (G_n) ( F(G_1, \ldots, G_n) \sim A(G_1, \ldots, G_n))$ [\P] is an axiom. This formula (V) has the purpose of replacing a formula which represents an individual predicate with a predicate variable (\cite{hilbert1938} p. 125).\footnote{Hilbert and Ackermann employ $\sim$ to denote the biconditional, 
$E F$ to denote existential second-order quantification, and $(G)$ to denote universal second-order quantification.   This passage thus expresses the full second-order comprehension scheme in the now familiar manner.}
\end{quote}

It is sometimes claimed that Supplement~IV of Volume 2 of \citep{hilbert1939} is the origin of full second-order Peano arithmetic (cf. \cite{simpson2009} p. 6) -- i.e. the system $\mathsf{Z}_2$ of Simpson's monograph.  However, the systems considered there contain neither the full nor a restriction version of the comprehension scheme. Rather, Hilbert and Bernays consider three second-order systems, called $H$, $K$, and $L$, respectively in sections A, F, G of Supplement IV. The system $L$ is in many ways similar to $\mathsf{Z}_2$, while the system $K$ differs from it in that the second-order entities have functions from natural numbers to natural numbers as their intended interpretation, instead of subsets of natural numbers. The system $H$ differs from system $K$ in that it uses a second-order variant of Hilbert's epsilon calculus, so that the first- and second-order quantifiers are defined in terms of the epsilon operator.\footnote{The epsilon operator and its associated calculus are of course less frequently studied today; but see \cite{leisenring1969a} for an overview.}

However, none of these systems have comprehension explicitly built-in; rather, what seems to go proxy for this are certain principles of explicit or stipulative definition (\cite{hilbert1939} pp. 454, 482, 490). They go proxy in the sense that, in their proof sketches of why the least upper bound principle is satisfied, they seem to be supposing that the stipulatively defined concepts fall within the range of the higher-order quantifiers.\footnote{More specifically, in reference to the proof-sketch on \cite{hilbert1939} pp. 463-464: the existence of a least upper bound of a non-empty bounded set $A$ of real numbers is explicitly presented in equation~(5) on the bottom of p. 463 in terms of a higher-order existential quantifier. In particular, the proof proceeds by presenting a ``definitional equation'' of a higher-order entity~$\nu$ at the bottom of p. 463. And the argument on p. 464 in equations (5a)-(5c) shows that $\nu$ is the least upper bound.} The situation in Supplement~IV is however less than clear because it is also suggested that the defined terms may be eliminated  (\cite{hilbert1939} pp. 455, 487). 

But of course second-order Peano arithmetic is not a conservative extension of this system with comprehension removed (or even restricted down to first-order comprehension). Rather, we now recognize that full comprehension decisively adds to the strength of the system and thus can by no means be regarded as a type of explicit definition, in any traditional sense of the term.  Moreover, a common gloss on the distinction between the systems studied in reverse mathematics is that they differ precisely in virtue of the ``set existence principles'' they contain -- a distinction which can be partially measured in terms of the inclusion of subschema of full comprehension.\footnote{But see \S\ref{sec:conclusions} for more on the qualification of ``partially''.}  In light of these complexities, it will be useful to say a bit more about how the understanding of comprehension principles evolved in relation to the development of second-order logic.   

Although \cite{hilbert1938} can reasonably be regarded as the first ``textbook'' treatment of second-order logic, the passage cited above from the second edition does not in fact represent the first statement of the second-order comprehension schema itself.  In fact, it is stated very clearly by G\"odel in his 1931 paper containing the Incompleteness Theorems, where the connection is that he proved the Incompleteness Theorem for a kind of simple type theory (\cite{godel1986} pp. 154-155, Axiom~IV, \cite{ferreiros1999} p. 355). Tarski also includes it in his 1935 paper on truth, but almost as an afterthought to his set-up of simple type theory. Tarski calls the comprehension schema a ``pseudodefinition,'' and he tells us
\begin{quote}
This term we owe to Le\'sniewski, who has drawn attention to the necessity of including pseudodefinitions among the axioms of the deductive sciences in those cases in which the formalization of the science does not admit the possibility of constructing suitable definitions. [\ldots] Pseudodefinitions can be regarded as a substitute for the \emph{axiom of reducibility} [\ldots] (\cite{tarski1956a} p. 223 fn, \cite{tarski1936a} p. 344 fn).
\end{quote}
\noindent This last line also occurs in G\"odel, who writes that the comprehension schema ``plays the role of the axiom of reducibility (the comprehension axiom of set theory)'' (\cite{godel1986} pp. 154-155).

Another intimation of the comprehension schema is provided by Ramsey in 1925. As with Hilbert and his collaborators, Ramsey suggested that the Axiom of Reducibility is dispensable, so long as the higher-order entities are conceived to be ``objective'' and do not depend ``on our methods of constructing them'' (\cite{ramsey1925a} p. 365, cf. \cite{parsons2002} pp. 381--382). To illustrate this, Ramsey considered a formula which we would write as $\varphi(x)\equiv \forall \; F \; \psi(F,x)$ and says that it determines a member of the range of second-order entities (\cite{ramsey1925a} p. 368).  Hence, as with Supplement~IV of \cite{hilbert1939}, while Ramsey does not explicitly state the full comprehension schema, he explicitly makes use of some of its immediate consequences.

Finally, it is perhaps worth mentioning an alternative formulation of comprehension which Church developed in the 1940s. In his \emph{Introduction to Mathematical Logic} from 1944, Church formalized comprehension in terms of the following substitution schema $(\forall \; F \; \Phi(F))\rightarrow \Phi(\psi(x)/F)$, where in the last term the expression $\Phi(\psi(x)/F)$ means ``substitute the formula $\psi(x)$ for the atomic formula $Fx$ in $\Phi$'' (cf. \cite{church1944a} Rule VIII$^{\prime}$ p. 100, \cite{church1956a} Rule 509 p. 297).\footnote{Of course, yet another alternative formalization comes in the way in which comprehension is handled in Church's simple theory \cite{church1940a}, where it is covered by $\lambda$-terms.} If one considers $\Phi(F)\equiv \exists \; G \; \forall x \; (Gx\leftrightarrow Fx)$, then the associated instance of the substitution schema straightforwardly implies the comprehension schema. And of course, by taking contrapositives, the schema is equivalent to $\Phi(\psi(x)/F) \rightarrow \exists \; F \; \Phi(F)$. Written in this way, the connection to the idea that the formula $\psi(x)$ determines a higher-order entity~$F$ becomes even more apparent.

But in its original non-contraposed formulation, the schema can end up looking a lot like the validity $(\forall \; x \; \varphi(x))\rightarrow \varphi(t)$ of predicate logic. Viewed from this perspective, the comprehension schema can take on the appearance of a tautology, which might explain why it took a comparatively long time for people to realize its import and strength.\footnote{This elision of substitution and comprehension is sometimes attributed to Frege (cf. \cite{demopoulos2005a} p. 131, \cite{sullivan2004a} p. 672).} Indeed, Henkin composed an entire paper as late as 1953 in which he noted that Church's substitution schema was equivalent to the comprehension schema (\cite{henkin1953a}). He thought that an advantage of formulating systems in terms of the comprehension schema was that it was natural to then consider ``certain subsystems'' obtained by weakening the comprehension schema, and for example one could entertain restrictions whose models were ``defined in a purely predicative way'' (\cite{henkin1953a} p. 207).

\subsection{Hilbert's finitism and primitive recursive arithmetic}\label{sec:hilfinPRA}

The significance of the axiom system $\mathsf{RCA}_0$ is often explained in relation to what (following \citet{kreisel1958}) is now called the \textit{Hilbert Program} -- i.e. the project by which Hilbert and his collaborators hoped to prove the consistency of analysis and portions of set theory using the mathematical resources made available by what they described as the \textit{finiten Standpunkt}.  But although this project can be seen as reaching its culmination in the \textsl{Grundlagen der Mathematik}, neither $\mathsf{RCA}_0$ nor any precise equivalent is described in this work.  To understand the connection between this system and the Hilbert program, it will thus be useful to consider a related system known as \textit{primitive recursive arithmetic} [$\mathsf{PRA}$], a version of which \textit{is} described in the \textit{Grundlagen}. 

Recall that $\mathsf{PRA}$ is the first-order theory whose language contains symbols for all primitive recursive functions and whose axioms consist of the defining equations for these functions together with the first-order induction scheme for quantifier-free formulas (cf. \cite{simpson2009} pp. 369-370).  There is a well-known proposal in the secondary literature on the Hilbert program originating with \citet{tait1968,tait1981} according to which finitary mathematics is characterized by the portion of mathematics which can be formalized within $\textsf{PRA}$.  And indeed the method of recursive function definition figures prominently not only in many of Hilbert and Bernays's expositions (e.g. \cite{hilbert2013a}, \cite{hilbert1926}, \cite{bernays1930}) which led up to \textsl{Grundlagen} but also with their description of what they call the \textit{finiten Standpunkt} in its first two chapters.   

The discussion of what Hilbert and Bernays call ``elementary number theory'' in chapter 2 of \citep{hilbert1934} characterizes its subject matter as finite sequences of symbols (i.e. \textit{numerals}) formed by the process of ``concretely terminating constructions'' (p. 21).  Primitive recursion is then explained as an ``abbreviated communication'' for the ``deconstruction of numerals'' and its role justified in terms of the fact that the relevant processes of deconstructing numerals into their parts may always be seen to terminate in a finite number of steps (pp. 26-27).  The application of induction and the least number principle  to \textit{decidable} formulas of elementary arithmetic is then given a similar justification in terms of finite search procedures (p. 23, pp. 34-35).  Hilbert and Bernays finally go on to summarize what they take to be characteristic of finitary mathematics is that it is confined to ``objects that are conceivable in principle'' and ``processes that can be effectively executed in principle''.  They thus conclude that ``it remains within the scope of a concrete treatment'' (p. 32).  

If such remarks are considered either in isolation or in regard to Hilbert and Bernays's earlier expositions, they might appear to lend credence to Tait's claim that they conceived of the \textit{finiten Standpunkt} as coinciding with what can be formalized in a system such as $\mathsf{PRA}$.  But if we are to view the \textit{Grundlagen} itself as one of the founding sources of reverse mathematics, then such a characterization becomes problematic for at least two reasons.  First, this work not only represents the final stage of the original development of the Hilbert Program, but it was also written at a time when Hilbert and Bernays were attempting to take into account not only G\"odel's completeness and incompleteness theorems (respectively in \cite{hilbert1939} \S4 and \S 5) but also the proof-theoretic work of Ackermann  and Gentzen.   Second, the characterization of the \textit{finiten Standpunkt} just recounted occurs at the beginning of the first volume of the \textit{Grundlagen}.  And in the course of the rest of its two volumes Hilbert and Bernays go on to consider a number of systems of formal arithmetic -- several of which properly extend $\mathsf{PRA}$ -- without ever explicitly stating that any of them captures their informal description of finitary mathematics.  

The first of these points has been the focus of recent discussion of Tait's claim.  For on the one hand, \citet{zach2003} observes that forms of transfinite recursion and recursion on higher types are used in Ackermann's dissertation (published as \cite{ackermann1925a}) and also that Hilbert appears to have accepted such methods as finitistic at this time.  And on the other, \citet{sieg2009} enumerates several instances in the \textit{Grundlagen} (most relating to formalization of consistency proofs) where Hilbert and Bernays employ methods which go beyond those which can be formalized in $\mathsf{PRA}$.  

In regard to the second point, it is notable that across the two volumes of the \textit{Grundlagen}, Hilbert and Bernays also consider at least eight distinct systems of formal first-order arithmetic.  These differ both as to their non-logical signature (the weakest contain only symbols for successor and less-than, while others contains primitive symbols for addition and multiplication), whether they contain schema allowing for additional of functions defined by primitive recursion, whether they contain an induction scheme, and if so, whether it is limited to quantifier-free formulas. Amongst these systems are (i)~the system $\mathsf{A}$ which is like Robinson's~$\mathsf{Q}$ in that it is induction-free and $\Sigma^0_1$-complete, but with the signature of just successor and less-than (\cite{hilbert1934} p. 263); (ii)~the $\mathsf{D}$ system like Presburger's arithmetic which has induction and defining equations for just addition (\cite{hilbert1934}, p. 357);  (iii)~the system $\mathsf{Z}$ which is identical to our contemporary axioms for $\mathsf{PA}$, minus any axioms for less-than (\cite{hilbert1934}, p. 371); (iv)~the system $\mathsf{Z}_{\mu}$ which extends $\mathsf{Z}$ with a form of the least number principle which Hilbert and Bernays show is equivalent to first-order induction (\cite{hilbert1939}, p. 293).

The transition between these systems is motivated both by both a desire to formalize greater fragments of number theory and in some cases to prove the consistency of the weaker systems and the stronger ones.  However Hilbert and Bernays observe that techniques which they have used to prove the consistency of the previous systems cannot be applied to $\mathsf{Z}$.  And after noting that each Diophantine equation may be expressed in $\mathsf{Z}$, they note that all recursive functions can be represented in this system, writing:
\begin{quote}
[\ldots] the formalism of system $\mathsf{Z}$ is not only, as we just found, in a position to formulate difficult problems of number theory, but it delivers more generally a formalization of the \emph{entirety of number theory}. Namely, in this formalism all functions are representable which can be introduced through recursive equations [\ldots] (\cite{hilbert1934} pp. 372-373).
\end{quote}

The system $\mathsf{Z}_{\mu}$ is introduced to aid in Hilbert and Bernay's exposition of G\"odel's completeness and incompleteness theorems in chapter 5 of \citeyearpar{hilbert1939}.  But although this ostensibly corresponds to the strongest theory considered prior to the introduction of the second-order theories in Supplement IV, it should finally be noted that Hilbert and Bernays explicitly demur from suggesting that even this system exhausts the resources of the finitary standpoint: 
\begin{quote}
The question arises as to whether finitary methods are in a position to exceed the domain of inferences formalizable in $\mathsf{Z}_{\mu}$. [\P] This question is admittedly, as so formulated, not precise; because we have introduced the expression ``finitary'' not as a sharply delimited endpoint, but rather as a designation of a methodological guideline, which would enable us to recognize certain kinds of concept formation and certain kinds of inferences as definitely finitary and others as definitely not finitary, but which however delivers no exact separating line between those which satisfy the demands of the finitary method and those which do not (\cite{hilbert1939} pp. 347-348, cf. \cite{sieg2009} p. 375).
\end{quote}

\subsection{$\mathsf{RCA}_0$ as a formal system}\label{sec:RCAformal}

Recall that $\mathsf{RCA}_0$ is derived from full second-order Peano arithmetic or $\mathsf{Z}_2$ by both restricting the scope of the comprehension scheme and also replacing the second-order induction axiom with the first-order induction scheme to $\Sigma^0_1$-formulas. Such a system would not have been directly formalizable in the framework of the \emph{Grundlagen} for, as we have just seen, comprehension was not distinguished as a separate logical principle there.  And although we have also seen that the full comprehension scheme was stated in Hilbert and Ackermann's textbook \cite{hilbert1938}, they do not consider the possibility of restricting it to specific subclasses of formulas.  Thus although subsequent work in proof theory attests to the close relationship of $\mathsf{RCA}_0$ to $\mathsf{PRA}$ -- and thus also to Tait's delineation of finitism -- the relevant results were not obtained until the 1970s.\footnote{In particular, by combining results of \citet{parsons1970} and \citet{friedman1976}, it is possible to show that $\mathsf{RCA}_0$ is conservative over the extension of $\mathsf{PRA}$ with first-order quantification theory (which is itself conservative over $\mathsf{PRA}$)  for $\Pi^0_2$-sentences. It hence follows that $\mathsf{RCA}_0$ is equiconsistent with $\mathsf{PRA}$ in the sense that any proof of a contradiction in the former could be transformed into a proof of a contradiction in the latter.  A series of well-known results about the system $\mathsf{WKL}_0$ considered in the next section suggest that similar points can be made about this system as well.  In particular, in the mid-1970s Harrington and Friedman showed $\mathsf{WKL}_0$ is $\Pi^1_1$-conservative over $\mathsf{RCA}_0$ and hence also $\Pi^0_2$ conservative over $\mathsf{PRA}$.  (Although Harrington and Friedman's model theoretic proofs of this result were never published, a related exposition was ultimately provided in Simpson's monograph \citep[IX.1-3]{simpson2009}.  Sieg obtained the same conservativity result for $\mathsf{PRA}$ by a proof theoretic argument \cite{sieg1985}.)  \citet{simpson1988a} subsequently suggested that these results contribute to a ``partial realization of Hilbert's program'' and also that $\mathsf{WKL}_0$ embodies a foundational standpoint which he labels \textsl{finitistic reductionism} \citep[p. 43]{simpson2009}.   As these claims have been widely discussed in the extant secondary literature on reverse mathematics -- e.g. \cite{feferman1988}, \cite{caldon2005}, \cite{burgess2010} -- we will not consider them further here. \label{note:conserv}}

A system equivalent to $\mathsf{RCA}_0$ was first introduced by Friedman in the abstracts \citeyearpar{friedman1976}.  But although his original formulation of the basic subsystems in \citep{friedman1975} was based on full first-order induction, no explanation for restricting induction is given in \citep{friedman1976}.  This topic is, however, discussed at length in \citep{friedman1983} wherein $\mathsf{RCA}_0$ appears to have first been formulated in its contemporary form.   One programmatic observation made there is that full induction is provably equivalent to the \textsl{bounded comprehension scheme} $\forall x \exists X(y \in X \leftrightarrow (y < x \ \wedge \ \phi(y)))$ over $\mathsf{RCA}_0$.  On this basis, Simpson subsequently remarked that since  ``the whole point of Reverse Mathematics is to prove ordinary mathematical theorems using only the weakest possible set~existence principle $\ldots$ the reverse mathematician is constrained to use full induction as sparingly as possible'' \citep[p. 150]{simpson1985a}.   It is additionally observed in both \citep{friedman1983} and \citep{simpson1985a} that the systems with restricted induction are more amenable to ordinal analysis, that their first-order parts of systems with restricted induction typically admit neater characterizations, and that while the use of restricted induction sometimes results in more complicated proofs, not only are few reversals to classical theorems lost in this manner, but their proofs sometimes yields improved quantifiers bounds.  

While these claims testify to the technical benefits of employing systems with restricted induction, they also provide part of the context of Simpson's subsequent discussion of $\mathsf{RCA}_0$ in relation to Bishop's constructive analysis.  Beginning with \citeyearpar{bishop1967}, Bishop presented a detailed constructive development of a large part of twentieth century analysis, inclusive of measure theory and the theory of Banach spaces.   His presentation avoids the use of formal systems and techniques from computability theory which he took to have hobbled the development of intuitionistic analysis after Brouwer.   It is clear, however, that Bishop's development is grounded in the supposition that analysis may be faithfully developed by treating natural numbers together with computable functions (or decidable sets) of numbers as basic objects.   This in turn inspired the formulation of several systems of constructive set theory (e.g. \citep{friedman1977}) wherein the Bishop's development of the basic number systems and their properties may be formalized.  

These papers in turn provided antecedents for the methods which are ultimately used to formalize analysis in $\mathsf{RCA}_0$  in \cite[\S II]{simpson2009}, but also provide the context for observation that the axioms of this system are ```constructive' in the sense that they are formally consistent with the statement that every total function from $\mathbb{N}$ into $\mathbb{N}$ is recursive'' \citep[p. 146]{friedman1983}.\footnote{This follows because it is easily seen that the $\omega$-model of $\mathsf{Z}_2$ whose second order domain consists of precisely the recursive sets satisfies $\mathsf{RCA}_0$ \citep[I.7.5]{simpson2009}.}  However, Simpson has also stressed (e.g. \citeyearpar[pp. 31-32]{simpson2009}) that the intention behind using  $\mathsf{RCA}_0$ as a base theory within Reverse Mathematics differs from Bishop's motivation of constructive analysis in both foundational and formal respects.   For on the one hand, the goal of calibrating the set~existence principles required to prove classical theorems is very different from the traditional constructivist goal of grounding mathematics in a theory of mental constructions and proving statements on the basis of the attendant understanding of the logical connectives.  And on the other, formal systems for constructive mathematics such as Heyting arithmetic (or its extension with higher types) are typically based on intuitionistic logic together with full first-order induction.  This is regarded as unproblematic both by \citet{bishop1967} and in classical expositions of intuitionism such as that of \citet{heyting1956}.  As we have seen, however,  $\mathsf{RCA}_0$ is based on classical logic together with a restricted induction scheme.   


\section{Weak K\"onig's Lemma and related systems}\label{sec:wkl_0}

The statement now often referred to as \textsl{K\"onig's Infinity Lemma} was first formulated by D\'enes K\"onig \citeyearpar[p. 122]{konig1927} as follows: ``If every point of a connected infinite graph has only finitely many edges going to it, then the graph contains an infinite path''.   K\"onig's isolation of this statement was the result of his attempts from 1914 onward to find a combinatorial proof of the Cantor-Schr\"oder-Bernstein Theorem.\footnote{K\"onig's father, Julius, had used a form of this result in his failed 1904 attempt to refute the Continuum Hypothesis by showing that the continuum could not be well-ordered and thus not equal to $\aleph_{\alpha}$ for any $\alpha$.   Although his proof was flawed, it led to a correct proof of what is now called \textsl{K\"onig's Theorem} in set theory -- i.e. if $\mathfrak{a}_i$ and $\mathfrak{b}_i$ are two families of cardinals indexed by $I$, such that $\mathfrak{a}_i < \mathfrak{b}_i$ for all $i \in I$ then $\Sigma_{i \in I} \mathfrak{a}_i < \Pi_{i \in I}\mathfrak{b}_i$.  Julius's original proof relied on a form of Cantor-Schr\"oder-Bernstein whose proof required the Axiom of Choice.  The work of  D\'enes which led to the Infinity Lemma was motivated by an attempt to determine the extent to which choice was necessary by reformulating the problem in graph-theoretic terms.  Although \citet[p. 171]{konig1990} observes that his proof of the Infinity Lemma still requires the Axiom of Choice when stated in full generality, he also notes this may be avoided in many of its applications.   See \citep{franchella1997} and \citep{hinkis2013} for more on K\"onig's use of the Lemma in his proof of the Cantor-Schr\"oder-Bernstein theorem.}  But by the time of his 1936 graph theory textbook, he had come to see the Infinity Lemma as a useful tool in its own right, writing that it has uses ``in the most diverse mathematical disciplines, since it often furnishes a useful method of carrying over certain results from the finite to the infinite'' (\cite{konig1990} p. 164, \cite{konig1936} p. 110).   Amongst these he lists the Heine-Borel Covering Lemma as well as a form of van der Waerden's Theorem, both of which he shows to be derivable from the Infinity Lemma.\footnote{The Covering Lemma states that the unit interval is compact with respect to the standard topology on the reals  -- i.e. ``Every covering of the closed unit interval $[0,1]$ by a sequence of open intervals has a finite subcovering''.   This statement  is now known to reverse to $\mathsf{WKL}_0$ over $\mathsf{RCA}_0$.   K\"onig's derivation of the Covering Lemma from the Infinity Lemma is similar to the proof given by \citet[IV.1.1]{simpson2009}.  He also states that this argument does not make use of the Bolzano-Weierstrass Theorem which notably reverses to the stronger system $\mathsf{ACA}_0$ over $\mathsf{RCA}_0$.  What K\"onig does not do, however, is to consider the possibility of proving a converse implication -- e.g. that the Infinity Lemma is itself derivable from Covering Lemma.} K\"onig also considered the restriction of the Infinity Lemma to \textsl{trees}  -- i.e. connected, acyclic graphs -- yielding the  familiar statement ``Every infinite, finitely branching tree has a infinite path'' [KL].

The statement now known as \textsl{Weak K\"onig's Lemma} [WKL] in reverse mathematics results from restricting KL not just to \textsl{binary} trees (i.e. trees in which each node has at most two children), but also to trees whose nodes are labeled only with the integers 0 and 1. The arithmetical formulation of this statement thus takes the following form: \textsl{every infinite subtree $T$ of the full binary tree $2^{< \mathbb{N}}$ has an infinite path}.\footnote{In this case, we assume that $T \subseteq \mathbb{N}$ via an encoding of finite sequences as natural numbers and that a \textsl{path} in $T$ is defined to be a function $f:\mathbb{N} \rightarrow \{0,1\}$ such that for all $k \in \mathbb{N}$, $\langle f(0),\ldots, f(k-1) \rangle \in T$.}   Although this may at first seem like an \textsl{ad hoc} restriction of the original principle, it is now known that WKL is sufficient to derive many statements of classical mathematics whose proofs, like that of the Covering Lemma, have traditionally been thought to require non-constructive choice principles.  However, this aspect of both KL and WKL originally came to light in the course of metamathematical investigations, to which we now turn. 
 
\subsection{The completeness theorem for classical predicate calculus}\label{sec:compCPC}

It is likely that G\"odel in his 1929 dissertation \citep{godel1986} was the first person to make metamathematical use of the Infinity Lemma in the course of his proof of the Completeness Theorem for the classical first-order predicate calculus.  Recall that G\"odel initially proved this result in the following form: if a first-order formula $\varphi$ is irrefutable, then $\varphi$ is satisfiable in some denumerable model.  

G\"odel's proof proceeds by constructing a sequence of finite models $\mathcal{M}_0, \mathcal{M}_1, \ldots$ which respectively satisfy formulas $\psi_0, \psi_1, \ldots$ obtained from the prenex normal form of $\varphi$.  These mimic the dependence of the bound variables of $\varphi$ in such a way that their joint satisfiability entails that of $\varphi$.  G\"odel showed that if $\varphi$ is satisfiable, then $\mathcal{M}_i \models \psi_i$ can be always be extended to $\mathcal{M}_{i+1} \models \psi_{i+1}$, but in only finitely many ways.  By treating these models as nodes in a tree determined by this extendability relation, an application of KL gives the existence of an infinite sequence of models $\mathcal{M}_0, \mathcal{M}_1, \ldots$  with the described properties.  A model $\mathcal{M} \models \varphi$ can now be obtained by letting the domain of $\mathcal{M}$ be a subset of the natural numbers and stipulating that a predicate $P(x_1,\ldots,x_k)$ appearing in $\varphi$ is satisfied by a $k$-tuple of natural numbers $n_{j_1},\ldots,n_{j_k}$ just in case there is an $i$ such that $\mathcal{M}_i \models P(n_{j_1},\ldots,n_{j_k})$.\footnote{G\"odel does not cite the Infinity Lemma by name in his proof, but rather says merely that the interpretation is obtained by ``familiar arguments''.  \citet[pp. 510-511]{heijenoort1967} reports that he was later told by G\"odel that these words were indeed intended to refer to the Infinity Lemma.  For a more detailed reconstruction of G\"odel's proof, see \citet[pp. 53-58]{godel1986}.}

In 1920 Skolem had previously described a similar construction in the course of his proof of the L\"owenheim-Skolem theorem in which he had employed the Axiom of Choice \cite[252--263]{heijenoort1967}. G\"odel observed this may be obtained as a corollary of his proof of completeness and similarly made no pretext that his proof was constructive.  He did, however, observe that the deductive completeness of the predicate calculus might be viewed as a form of decidability in the sense that it demonstrates that ``every expression $\ldots$ either can be recognized as valid through finitely many inferences or its validity can be refuted by a counterexample'' \citeyearpar[p. 63]{godel1986}.  But he also observed that in order for such a result to bear on the completeness of intuitionistic logic, the assertion that a countermodel exists would itself need to be proven constructively.

The question of whether the completeness theorem admits a constructive proof was considered more explicitly by Hilbert and Bernays \citeyearpar{hilbert1939}.   As we have seen, G\"odel had already observed that an irrefutable formula $\varphi$ possesses an \textsl{arithmetical model} -- i.e. one whose domain consists of a subset of the natural numbers, and whose predicate and function symbols are interpreted as sets of natural numbers of the appropriate arities. Examination of G\"odel's proof also makes clear that the arithmetical formulas $P^*_1,\ldots, P^*_n$ can be constructed uniformly from the formula $\varphi$.   In light of this, \citet[pp. 189-190]{hilbert1939} went on to introduce the notion of an \textsl{effectively satisfiable} [\textsl{effektiv Erf\"ullbar}] formula -- i.e. one which upon being put into prenex normal form can be transformed by effectively replacing atomic formulas with truth values and formulas containing free variables with computable [\textsl{berechenbar}] number theoretic predicates so that each substitution instance with numerals is made true in the standard model.   They thought that such an interpretation would constitute  a ``finite sharpening'' (p. \citeyearpar{hilbert1939}, p. 191) of the Completeness Theorem.  But upon introducing the definition of effective satisfiability, they go on to conjecture that completeness would \textsl{fail} if this notion were to be substituted for the traditional (non-effective) definition of satisfiability in its statement.

In light of Church and Turing's work on the \textsl{Entscheidungsproblem} (which Hilbert and Bernays take into account in Supplement II of \citeyearpar{hilbert1939}) such a conjecture would certainly have been reasonable by the late 1930s.   But it was not fully resolved until the early 1950s in virtue of work which is now thought of as contributing more directly to computability theory than it is to model theory.  In particular, it appears to have been \citet[p. 268]{kreisel1950} who first explicitly observed that a special case of the Infinity Lemma is sufficient for the proof of the Completeness Theorem.  In the same paper he also uses Hilbert and Bernays's arithmetization of G\"odel's proof to construct a statement which is formally independent of a subsystem $\mathsf{S}$ of G\"odel Bernays set theory (considered in \S \ref{sec:aca0formal}) above and suggested on this basis that in no arithmetical model of this theory could the membership relation $\in$ receive a recursive interpretation.\footnote{As Wang observed in his review of \citep{kreisel1950},  Kreisel's proof doesn't actually yield this result but (in effect) only the weaker statement that the definition of $\in^*$ produced by applying the Arithmetized Completeness Theorem to $\mathsf{GB}$ can never be $\Delta^0_1$.  This anticipates the later work (summarized in the introductory remarks to \citep{godel1986}) which collectively showed that any consistent statement has an arithmetical model in which all of its relations are $\Delta^0_2$-definable and also that this statement fails for $\Sigma^0_1 \cup \Pi^0_1$.} 

\citet{kreisel1953a} also formulated the relevant principle as a statement which can be expressed in a second-order extension of Hilbert and Bernay's system $\mathsf{Z}_{\mu}$.\footnote{The statement is of the following form: any infinite, finitely branching tree consisting of finite sequences of integers $\leq k$ (for some fixed $k \in \mathbb{N}$) which is also determined by a primitive recursive predicate of finite sequences has a path which is given by a term $t(x)$ of the form $\mu y.\forall z R(x,y,z)$ where $R(x,y,z)$ is itself a primitive recursive predicate \cite[p. 125-126]{kreisel1953a}.}  He did not, however, carry out the formalization of G\"odel's original argument from this principle.  However, the fact that it suffices to consider only subtrees of $2^{< \mathbb{N}}$ is evident from Kleene's \citeyearpar{kleene1952} reformulation of G\"odel's proof based on the method of maximally consistent sets introduced by \citet{henkin1949}.  

Recall that in this construction we assume an enumeration of all sentences $\psi_1, \psi_2, \ldots$ in the language of $\varphi$ augmented with with new constants $c_i$ as well as axioms $\exists x \psi_i(x) \rightarrow \psi(c_i)$.  Letting $\chi_0 = \varphi$ we then consider at the $(i+1)$-st stage the result of successfully setting $\chi_{i+1} = \psi_{i+1}$ if $\{\chi_0,\ldots,\chi_{i}\} \cup \{\psi_{i+1}\}$ is consistent and $\chi_{i+1} = \neg \psi_{i+1}$ otherwise.   Relative to the given enumeration, finite sequences of this form can be represented as finite binary sequences $\sigma \in 2^{< \mathbb{N}}$ where $\sigma(i) = 0$ if $\chi_i = \psi$ and $\sigma(i) = 1$ if $\chi_i = \neg \psi_i$.  Suppose we consider at the $i$th stage all sequences satisfying the following predicate:
\begin{eqnarray*}
S_{\varphi}(\sigma) & = & \text{there is no proof in the predicate calculus of a contradiction } \\
& & \text{ of length less than that of $\sigma$ from the set encoded by $\sigma$}
\end{eqnarray*}
It is evident that the set of sequences satisfying $S(\sigma)$ is a tree $T_{\varphi} \subseteq 2^{< \mathbb{N}}$ and that $T_{\varphi}$ is infinite just in case $\varphi$ is irrefutable.  Moreover, from an infinite path through $T_{\varphi}$ a model of $\varphi$ can be constructed in the familiar manner of the Henkin construction.  

As finite binary sequences can be encoded as natural numbers and the quantifier over proofs in the definition of $S_{\varphi}(\sigma)$ is bounded, it follows that $T_{\varphi}$ can be defined by a primitive recursive predicate.  Generalizing one step further, let $\Gamma$ be an arbitrary recursively axiomatized theory.  We may consider the predicate $S_{\Gamma}(\sigma)$ defined analogously $S_{\varphi}(\sigma)$ where we consider proofs of length less than $|\sigma|$ from the axioms of $\Gamma$.  In this case $T_{\Gamma}$ is infinite just in case $\Gamma$ is consistent and is definable arithmetically by a recursive (i.e. $\Delta^0_1$) predicate in the language of first-order arithmetic.  Consideration of recursive trees of this form led to another development in mathematical logic which anticipated the isolation of $\mathsf{WKL}_0$ as a formal system -- i.e. the formulation of the so-called \textsl{basis theorems} in computability theory.

\subsection{The constructive failure of K\"onig's Lemma and the basis theorems}\label{sec:basisthm}

In the terminology of contemporary computability theory, a class of sets $\mathcal{S}$ of natural numbers is known as a $\Pi^0_1$-\textsl{class} if it is definable by some $\Pi^0_1$-formula. There is a well-known representation theorem, provable in $\mathsf{RCA}_0$, that every $\Pi^0_1$-class is representable as the set of paths through some recursive binary tree $T \subseteq 2^{< \mathbb{N}}$. WKL hence expresses the non-emptiness of $\Pi^0_1$-classes whose underlying tree is infinite.  By applying the construction just described,  \citet{kleene1952a} showed that if $\Gamma$ is a recursively axiomatizable essentially undecidable theory (such as first-order Peano arithmetic), then $T_{\Gamma}$ can contain no infinite recursive path.  For in this case it can be shown that the sets $A$ and $B$ consisting of the G\"odel numbers of sentences provable and refutable in $\Gamma$ form a pair of computably inseparable sets.  On the other hand, any path $f:\mathbb{N} \rightarrow \{0,1\}$ through $T_{\Gamma}$ corresponds to the characteristic function of a set separating $A$ from $B$ -- i.e. if $C = \{n \ : \ f(n) = 1\}$ then $C \supseteq A$ and $C \cap B = \emptyset$.  Thus $f(x)$ cannot be recursive. This result is interesting in its own right, as it is a failure of a \emph{separation principle}, to which we shall return in \S\ref{sec:addsteel}.

For present purposes, however, the import of this result is that there exist non-empty $\Pi^0_1$-classes with no recursive members. \citet{kleene1952a} originally obtained this result in the context of investigating whether Brouwer's \textsl{Fan Theorem} is consistent with the assumption that all choice sequences are recursive.  In the relevant case, the Fan Theorem asserts that if every tree $T \subseteq 2^{< \mathbb{N}}$ is such that each path $f(x)$ through $T$ (notation: $f \in [T]$) has an initial segment which satisfies some property $A$, then there is a uniform bound on the length of the initial segments of $f(x)$ (notation: $f \restriction x$) at which this property is satisfied:
\[\forall \; f \; \in \; [T] \; \exists \; x \; \varphi(f \restriction x) \rightarrow \exists \; y \; \forall f \; \in [T] \; \exists \; x \leq y \; \varphi(f \restriction x)\] 
This principle expresses a form of compactness of the intuitionistic continuum which Brouwer sought to preserve in his development of intuitionistic analysis, such as his proof that every continuous real-valued function on a closed interval is uniformly continuous \citep{brouwer1927}.  

It is readily seen that the Fan Theorem is classically equivalent to K\"onig's Lemma.  But Kleene's result has also led to some theorists working within the intuitionistic and constructive traditions to regard the Fan Theorem with suspicion.  In particular, \citet{beth1947,beth1956} had made essential use of the Fan Theorem in his completeness proof for intuitionistic first-order. In conjunction with G\"odel and Dyson, Kreisel \citep{dyson1961,kreisel1962,kreisel1970} then established a series of results showing that the formalized statement of completeness for the intuitionistic predicate calculus with respect to what are now called \textit{Beth models} entails Markov's Principle and the negation of the intuitionistic  Church's Thesis (one form of which states  ``every function $f:\mathbb{N} \rightarrow \mathbb{N}$ is recursive'').  On this basis the constructive significance of formal completeness proofs has also thereby repeatedly been called into question -- cf. \citep{dalen1973, troelstra1988a}.\footnote{At the same time, a minority view takes the moral of these equivalences to be that the intuitionist ``has to admit the possibility of infinite sequences of natural numbers that may be effectively calculated but are not given by an algorithm'' \citep[p. 635]{veldman2014}.}

By further examination of the Hilbert and Bernays's arithmetization of the G\"odel completeness proof, \citet[\S 72]{kleene1952} and \citet{hasenjaeger1953} also showed that if $\Gamma$ is recursively axiomatizable, then $T_{\Gamma}$ must contain a path which is $\Delta^0_2$-definable in the language of first-order arithmetic.  As every infinite recursive binary tree can be represented in this form for an appropriate theory $\Gamma$, the class of sets defined by $\Delta^0_2$-predicates is said to serve as a \textsl{basis} for non-empty $\Pi^0_1$-classes -- i.e. every non-empty $\Pi^0_1$-class has a member among the class of sets defined by $\Delta^0_2$-predicates. \citet{kleene1959} attributed the following general definition of a basis to Kreisel:  a class of $\mathcal{C}$ of subsets of $\mathbb{N}$ is a basis for a class $\Theta$ of second-order formulas containing $X$ free just in case if each instantiated predicate $\varphi(X) \in \Theta$ is such that $\exists X[X \in \mathcal{C} \ \wedge \varphi(X)]$ (see also \citep[\S 7.11]{shoenfield1967}). He additionally writes that the significance of the fact that $\mathcal{C}$ is a basis for $\Theta$ is that definitions given by formulas of this form ``mean the same to persons with various universes of [sets], so long as each person's universe includes at least $\mathcal{C}$'' \citeyearpar[p. 24]{kleene1959}.  

By a well-known result of \cite{post1948}, a set $A \subseteq \mathbb{N}$ is $\Delta^0_2$-definable just in case $A$ is Turing reducible to the halting set $K$ (i.e. $A \leq_T K$) or, equivalently, that the Turing degree of $A$ is less than or equal to $\emptyset'$ (i.e. $\textrm{deg}(A) \leq_T \emptyset'$).  It thus follows that the class of sets of Turing degree $\leq_T \emptyset'$ form a basis for non-empty $\Pi^0_1$-classes (the \textsl{Kreisel Basis Theorem}).  This result was subsequently strengthened by \citet{shoenfield1960} to show that the class of sets of Turing degree $<_T \emptyset'$ forms a basis for non-empty $\Pi^0_1$-classes (the \textsl{Shoenfield basis theorem}) and again by \citet{jockusch1972}  to shows that the class of sets $A$ such that $\textrm{deg}(A') \leq_T \textrm{0}'$ form a basis for non-empty $\Pi^0_1$-classes (the \textsl{low basis theorem}).  

These results anticipate the formulation of $\textsf{WKL}_0$ as a formal system in the sense that they provide natural computability-theoretic characterizations of its $\omega$-models.  In particular, recall that a non-empty set $\mathcal{S} \subseteq 2^{\mathbb{N}}$ is a \textsl{Turing ideal} just in case it is closed under effective join (i.e. if $A,B \in \mathcal{S}$, then $A \oplus B \in \mathcal{S}$) and Turing reducibility (i.e. if $A \in \mathcal{S}$ and $B \leq_T A$, then $B \in \mathcal{S}$) and that if $\mathcal{S}$ is a Turing ideal, then $\langle \mathbb{N}, \mathcal{S},+,\cdot,0,1,< \rangle$ is model of $\mathsf{RCA}_0$.  If $\mathcal{S}$ is additionally closed under the condition \textsl{if $A \in \mathcal{S}$ codes an infinite subtree of $2^{< \mathbb{N}}$, then there exists $B \in \mathcal{S}$ which is a path through $A$} then $\mathcal{S}$ is known as a \textsl{Scott set}.  It is easy to see that $\omega$-models of the language of second-order arithmetic whose second-order domains $\mathcal{S}$ are Scott sets are models of $\mathsf{WKL}_0$.  The results just summarized thus provide examples of such models where $\mathcal{S}$ is obtained as the Turing ideal generated by some basis for a non-empty $\Pi^0_1$-classes.  

Another source of models of $\mathsf{WKL}_0$ which anticipated the formulation of  $\textsf{WKL}_0$ is provided by Scott's \citeyearpar{scott1962} work on classes of sets bi-enumerable in complete consistent extensions of first-order Peano arithmetic.  Such classes may be characterized semantically as the so-called \textit{standard systems} of nonstandard models $M$ $\mathsf{PA}$ -- i.e. $\text{SSy}(M) = \{A \subseteq \mathbb{N} \ : \exists B \in \textrm{Def}(M)[A = B \ \cap \ \mathbb{N}]\}$ where $\text{Def}(M)$ denotes the set of sets which are definable with parameters in $M$.  Friedman observed that the standard system of any countable nonstandard model of $\mathsf{PA}$ is a Scott set and also that any Scott set is the standard system of some such model \citep[p. 541-542]{friedman1973}.  It thus follows from another observation of \citet[p. 238]{friedman1975} that the $\omega$-models of $\textsf{WKL}_0$ are precisely those whose second-order part is $\text{SSy}(M)$ for some countable nonstandard $M \models \mathsf{PA}$.

\subsection{$\mathsf{WKL}_0$ as a formal system}\label{sec:WKLformal}

Recall that the axiomatic theory $\mathsf{WKL}_0$ consists of $\mathsf{RCA}_0$ together with the arithmetical formulation of WKL described above.  Although this theory was first introduced by \citet{friedman1975}, such a system is described informally by Kreisel, Mints, and Simpson \citeyearpar{kreisel1975}.  This paper is devoted to assessing the extent to which  ``abstract language'' pertaining to the existence of infinite sets is necessary for either stating number theoretic results or for providing comprehensible proofs of arithmetical theorems whose statements themselves do not require such language.  The authors also explicitly discuss the methodology of using subsystems of second-order arithmetic to study the set~existence principle implicit in different foundational standpoints:
\begin{quote}
[T]here is a logical view which requires restrictions, for example because of (genuine or ethereal) doubts about the existence of sets having certain formal properties $\ldots$ It certainly can do no harm to have some idea of the consequence of a given `view', for example, for the class of, say, number theoretic theorems provable by means of [an] abstract principle. \citep[p. 116]{kreisel1975} 
\end{quote}

It is also in this paper that it was first proposed that KL and WKL could be understood as \textsl{axioms} which can be added to weak arithmetical system so as to formalize ``abstract'' statements about infinite sets.  The authors explicitly distinguish between KL, the restriction of KL to binary trees labeled with integers of unbounded size (KL$^-$), and WKL.  Upon first observing that there is no reason to suspect that that WKL is as strong as KL$^-$, they then go on to identify the $\Delta^0_1$-comprehension schema and describe a system similar to $\mathsf{RCA}_0$ as a potential base theory for investigating such claims.  They record that Friedman had already shown that the result of adding KL$^-$ to $\Delta^0_1$ plus ``closure under a few primitive recursive operations'' yields the full arithmetical comprehension principle (i.e. presumably a system coincident with $\mathsf{ACA}_0$) and is thus not conservative over this system.  However, Kreisel, Mints, and Simpson also observe (pp. 124-125) that the result of adding KL and  $\Delta^0_1$-CA to $\mathsf{PA}$ yields a conservative extension.\footnote{The equivalence of both KL and KL$^-$ to arithmetical comprehension over $\mathsf{RCA}_0$ is presented as Theorem III.7.2 in \cite{simpson2009} and attributed to \citet{friedman1975}.  The conservativity result then follows since $\mathsf{ACA}_0$ is a conservative extension of $\mathsf{PA}$.}  

What Kreisel, Mints, and Simpson do not do, however, is to conjecture the Friedman-Harrington result that $\textsf{WKL}_0$ is conservative over $\textsf{RCA}_0$ for $\Pi^1_1$-formulas (see note \ref{note:conserv} above).  And although this result was also stated without proof by \citet{friedman1975}, it is not this specific feature of $\textsf{WKL}_0$ which he takes to illustrate why this is a natural subsystem to consider in the development of reverse mathematics. Rather he suggests that $\textsf{WKL}_0$ is an example of the following general theme: ``Much more is needed to define explicitly a hard-to-define set of integers than merely to prove [its] existence''  \citeyearpar[p. 235]{friedman1975}.  

The computability-theoretic aspects of $\textsf{WKL}_0$ discussed above illustrate why this is so.  For as Friedman observes, $\textsf{WKL}_0$ is sufficient to prove that a non-computable set exists (as follows from the fact that we can formalize Kleene's argument in this system). But this system is not strong enough to prove the existence of the Turing jump of an arbitrary set (as follows from the fact that $\mathsf{WKL}_0$ has an $\omega$-model consisting of just low sets).  However, since $X'$ is $\Sigma^0_1$-definable relative to $X$, an application of arithmetical comprehension yields that $\mathsf{ACA}_0$ proves that $X'$ always exists whenever $X$ does.


\section{Arithmetical transfinite recursion and countable ordinals}\label{sec:ATR0}


The system $\mathsf{ATR}_0$ formalizes, and reverses to, certain elements of classical descriptive set theory. Appropriately, the intellectual origins of this system lie in the early history of descriptive set theory in Borel and Lusin, which we set out in \S\ref{subsec:borel}. Then, in \S\ref{sec:addsteel}-\S\ref{sec:ATRformal}, we describe the more immediate antecedents to Friedman's axiomatization of $\mathsf{ATR}_0$: namely, the effectivization of descriptive set theory by Addison and Kreisel, and the development of hierarchies by Harrison.

\subsection{Borel, Lusin, and countable ordinals}\label{subsec:borel}

A salient aspect of the early history of descriptive set theory was a skepticism about the existence of ordinals, in spite of their ostensible presence in certain core concepts like the hierarchy of Borel sets. For instance, the title of Souslin's famous paper which contains his result that the Borel sets are precisely those sets which are both analytic and co-analytic, is ``On a definition of the Borel measurable sets without transfinite numbers'' (\cite{souslin1917a}). In this same spirit, Kuratowski published a paper indicating how to avoid transfinite ordinals in certain constructions, saying that ``in reasoning with transfinite numbers one makes implicit use of their existence; now the reduction of the system employed in the demonstration is desirable from the point of logic and mathematics'' (\cite{kuratowski1922a} p. 77). Kuratowski's method is similar to how we might alternatively define the Borel sets to be the smallest collection satisfying certain closure properties.

Now, as for its sources, one source of skepticism about ordinals was due to Borel, who prior to the paradoxes in 1898 argued that Cantor's second principle of generation of ordinals, namely the taking of supremums, could not of itself generate an uncountable ordinal. Retrospectively, of course, we can see that Borel had a point: for if one formalizes the second principle as saying that the supremum of any set of ordinals exists, then this principle would be validated in models like $V_{\omega_1}$, wherein all ordinals are countable. The way that Borel put the point was that ``the second principle of formation could only make us acquire cardinals which we already have'' and he went onto add: ``and it seems doubtful that we have an idea sufficiently precise of what could be a cardinal exceeding the countable'' (\cite{borel1898} p. 122, \cite{borel1914} p. 122, cf. \cite{gispert1995a} p. 61).  But Borel's conclusion is much more agnostic than that of Brouwer, who in his 1907 dissertation gave exactly the same argument but with the conclusion that ``Cantor's second number class does not exist,'': ``[\ldots] it cannot be thought of'' and ``[\ldots] it cannot be mathematically constructed'' (\cite{brouwer1975a} p. 81).\footnote{See \cite{troelstra1982a} \S\S{2-4} for more on the influence of Borel on Brouwer.}

Another source of skepticism about ordinals was the Burali-Forti paradox. While Hadamard thought that this paradox was no different than the paradoxes which initially beguiled other fruitful mathematical concepts like the negative and complex numbers (\cite{hadamard1905b} p. 242, cf. \cite{garciadiego1992a} p. 139), Poincar\'e used the paradox to inveigh against mathematics which was not sufficiently rooted in intuition (\cite{poincare1905a} \S{8} pp. 824-825, \cite{poincare1906a} \S{7} pp. 303-305, cf. \cite{moore1981a} pp. 340-342). Lusin then went onto use Burali-Forti to argue against the claim that we have an intuition of ordinal numbers on the basis of the familiar representation of small countable ordinals (as various products and sums of $\omega$ and finite ordinals):
\begin{quote}
In effect, if we make (or believe to make) an image perfectly clear of the totality of countable ordinals [\P] $0$, $1$, $2$, \ldots, $\omega$, $\omega+1$, \ldots, $\alpha$, [\P] we see with the same clarity the totality of \emph{all} the transfinite numbers, and by the reasoning of Burali-Forti we see that this totality is logically contradictory in itself (\cite{lusin1930a} p. 26). 
 \end{quote}
\noindent Lusin and Borel then suggested viewing ordinals merely as an ``abbreviated notation'' for ``the order in which must be effected a countable number of operations'' (\cite{lusin1930a} p. 29, cf. \cite{borel1914} p. 231).

On Lusin's view, this explicit construction is crucial to Borel sets: ``It is the order of the intermediary sets which is the veritable nerve of the \emph{constructive} definition of Borel measurable sets'' (\cite{lusin1930a} p. 29). By contrast, non-Borel projective sets have a kind of secondary status, perhaps similar to how we might say that something is merely a class (and not a set) or merely given by a formula (as opposed to determining a second-order object). Speaking of projective non-Borel sets, Lusin writes that he
\begin{quote}
[\ldots] adopts the empiricist point of view and is inclined to consider the examples constructed by him as forms of words and not as defining objects genuinely completed, but only virtual objects (p. 322).
\end{quote}
\noindent This view in 1930 is far more tempered than the view espoused by Lusin in the original 1925 papers, in which many of the results from the 1930 book were first obtained. In particular, due to his skepticism about the totality of countable ordinals (and referencing Baire's famous statement from the \emph{cinq lettres} (cf. \cite{hadamard1905a} p. 264)), Lusin says in 1925 that 
\begin{quote}
Baire wrote, 20 years ago, in his letter to Hadamard, ``of a set which is given, it is in my view false to consider the subsets of this set as given.'' It seems to us that this assertion could be extended in the following manner: a set being given, it is false to consider its \emph{complement} as given (\cite{lusin1925e} p. 281, italics added).
\end{quote}
\noindent Lusin is here referring to the now familiar decomposition of a \emph{co}analytic set $X$ as $X=\bigcup_{\alpha< \beta} X_{\alpha}$, where the $X_{\alpha}$ are pairwise disjoint Borel and $\beta\leq \omega_1$ (cf. \cite{lusin1930a} pp. 204-205; cf. \cite{kechris1995a} p. 269, \cite{moschovakis2009a} p. 162). Thus in~1925, Lusin had voiced the claim that set~existence is not necessarily closed under \emph{complementation}, in part due to his skepticism about the totality of all countable ordinals.

By contrast, in 1930 Lusin suggested that there was a possibility of a kind of experimental confirmation of the totality of all countable ordinals. In particular, Lusin formulated the problem of determining whether every coanalytic set is countable or of the cardinality of the continuum \citep[p. 295, Problem 1]{lusin1930a}. He thought that the interest in this lies in a \emph{negative} resolution, of the form: there exists an uncountable coanalytic set without a non-empty perfect subset-- in the contemporary parlance, an uncountable thin coanalytic set (cf. \cite{moschovakis2009a} pp. 187, 212).\footnote{Of course, strong forms of determinacy axioms imply the \emph{positive} resolution. See \cite{jech2003} p. 629.} For, as mentioned above, Lusin was able to write a coanalytic set $X$ as $X=\bigcup_{\alpha< \beta} X_{\alpha}$, where the $X_{\alpha}$ are pairwise disjoint Borel and $\beta\leq \omega_1$. And if $X$ did not include a non-empty perfect subset, then neither would any of the $X_{\alpha}$, and hence by the perfect set theorem for Borel sets, each $X_{\alpha}$ would be countable, which under the hypothesis that $X$ was uncountable would entail that $\beta=\omega_1$. That is, this circumstance would result in the least uncountable ordinal~$\omega_1$ being rather concretely realized in terms of a partition of a coanalytic set of reals into $\omega_1$-many disjoint countable Borel parts. Lusin puts the significance of this as follows:
\begin{quote}
Therefore, \emph{if this case, which is logically possible, is practically real, one could affirm that the existence of all countable ordinals is an empirical fact}. (\cite{lusin1930a} p. 295).
\end{quote}

Lusin articulated a method for solving this problem, which he called \emph{the method of resolvants}. As he notes, its origin was in a remark from Borel's 1908 lecture at the International Congress of Mathematicians, wherein Borel sets out his view of the arithmetic continuum (which he opposes to the \emph{geometric continuum}, which is given in intuition):
\begin{quote}
The continuum never appears given in its totality, from the arithmetical point of view; each of its elements could be defined (or at least, there are none of its elements of which we could actually affirm that they could not be defined) [[Footnote]: Here is, to give an idea of my point of view, a problem which appears to be the most important in the arithmetic theory of the continuum: is it or not possible to define a set $E$ such that one could not name any individual element of this set $E$, that is to say, such that one could not distinguish without ambiguity this from all the other elements of $E$? (\cite{borel1909} p. 17, \cite{borel1914} pp. 161-162).]
\end{quote}
\noindent Hence Borel raises the question of whether definable non-empty sets necessarily have definable members. 

In Lusin's hands, this idea became identified with a method of solving problems. As Lusin's student Keldysh put it, Lusin ``says that a problem has been `\emph{mis en r\'esolvante}' if there exists a point set $E$ such that the problem can be resolved affirmatively if a point in the set $E$ can be specified, and resolved negatively if it can be proved that $E$ is empty'' (\cite{keldysh1974} p. 185, cf. \cite{kanovei2003} p. 866, \cite{lusin1930a} p. 293). Lusin then showed that the resolvant of the problem of the existence of uncountable thin coanalytic sets is a projective set (cf. \cite{lusin1930a} p. 295, \cite{keldysh1974} p. 185). Part of Lusin's prediction that the famous problems of descriptive set theory would remain unsolvable was that ``one can neither name an individual point in such sets [projective sets], nor know if there `exist' points in such a set, nor know their properties'' (\cite{lusin1930a} p. 303). That is, while the projective sets are definable, Lusin predicted that they do not all have definable elements, and hence are precisely the kind of sets which Borel's ``most important'' problem asked after.

\subsection{Effectivizing descriptive set theory}\label{sec:addsteel}

These ideas of Borel and Lusin found their way into latter-day developments primarily through the work of Addison. His dissertation, written under Kleene and finished in 1954, is in part occupied with formally defining the basic concepts of what is now commonly referred to as effective descriptive set theory. And at the outset of his dissertation, Addison writes that \citet{lusin1930a} ``has been our constant companion and guide'' (\cite{addison1954a} p. 4).

Addison points out that Borel had been one of the first to suggest the idea of some kind of effectivization of descriptive set theory, or at least the Borel sets. But writing in 1914, Borel did so not in terms of computation, but in a less formal manner. For instance, he calls a real number~$r$ \emph{calcuable} if given an $n>0$ one ``knows'' a rational~$q$ such that $\left|r-q\right|<\frac{1}{n}$. Amplifying upon this in a footnote, he tells us that what's essential in such knowledge is that ``each of the operations [in the calculation] must be executable in a finite amount of time, by a method certain and without ambiguity'' (\cite{borel1914} p. 219 fn). He extends this to functions from reals to reals, saying that such a function~$f$ is calculable when~$f(r)$ is calculable for all calculable~$r$ (\cite{borel1914} pp. 223--224). He then says that a set is \emph{bien d\'efinis} (literally: well-defined), if its characteristic function is calculable (\cite{borel1914} p. 225), and that sets \emph{bien d\'efini} are ``precisely'' those which in the first 1898 edition of his text he called \emph{measurable}, and which Lebesgue later renamed \emph{B-measurable}, where of course the ``B'' stands for ``Borel'' (\cite{borel1914} p. 226).\footnote{The first 1898 edition of Borel's text \citep{borel1898} contains no discussion of calculability, and so one should take Borel's ``precisely'' with a grain of salt. The discussion of calculability first occurs in the 1914 edition of Borel's text in the context of Richard's paradox (cf. \cite{borel1914} pp. 162 ff).}

While prescient, there is obviously much that is informal in this aspect of Borel's work, and Addison took his task to be to formalize Borel's ideas using tools from computability theory (\cite{addison1954a} pp. 43 ff). He defines effective Borel subsets of Baire space, along with natural number indices of them, by transfinite recursion, defining a sequence $K_{\alpha}=K_{\alpha}^+\cup K_{\alpha}^-$ in analogue to the usual definition of sequence $B_{\alpha}$ of the Borel subsets. In particular, the clopen basis $K_0=K_0^+=K_0^-$ of Baire space is given by the functions that pass through a given finite sequence, and Addison presents an effective coding of these which serve as the indices. Then a set said to be in $K_{\alpha}^+$ (resp. $K_{\alpha}^-$) with index $e$ if $e$ is an index for a total computable function such that the set can be written as a union  (resp. intersection) of sets $X_n$ such that $X_n$ has index $\{e\}(n)$, where as usual this denotes the action of the $e$-th program on input~$n$. Thus $K_1^+$ is what we now recognize as the effective open sets and $K_1^-$ as the effective closed sets, or what we often refer to simply as the $\Pi^0_1$-classes (cf. \S\ref{sec:basisthm}). Finally, he notes that whereas the classical hierarchy of Borel sets is equal to $\bigcup_{\alpha<\omega_1} B_{\alpha}$, the effective hierarchy is equal to $\bigcup_{\alpha<\omega_1^{CK}} K_{\alpha}$, where $\omega_1^{CK}$ denotes the least non-computable ordinal (cf. \cite{addison1954a} pp. 46 ff).

This general construction from Addison's dissertation was never published. However, the classical treatment of effective descriptive set theory is still Chapter 3 of Moschovakis' descriptive set theory textbook, where we are told that Addison ``[\ldots] initiated the development of the unified treatment we are presenting here'' (\cite{moschovakis1980a} p. 118, \cite{moschovakis2009a} p. 88). However, all of Addison's subsequent research was animated by a related aspect of the dissertation, where he explores the analogies between computability theory and descriptive set theory, particularly as regards their separation principles (cf. \cite{addison1959b}, \cite{addison1962a} \S{3}, \cite{addison1968a}, \cite{addison2004a} \S{3}).

Recall that a class of sets is said to satisfy the \emph{separation property} if for any two disjoint members~$A$, $B$ from the class there is $C$ such that $A\subseteq C$ and $B\cap C=\emptyset$ and such that both $C$ and its complement belong to the class in question. In \S\ref{sec:basisthm} we had occasion to note Kleene's result that the $\Sigma^0_1$-definable sets of natural numbers do \emph{not} satisfy the separation property. This contrasts to Lusin's famous result that the analytic sets do satisfy the separation property (\cite{lusin1930a} p. 156, \cite{moschovakis2009a} p. 156). Since both the analytic sets and the $\Sigma^0_1$-definable sets of natural numbers are defined in terms of an existential quantifier, this initially seemed like it suggested that the two hierarchies were very dissimilar.

However, Addison's dissertation culminates in an explanation of why the analogy should rather be: analytic sets correspond to $\Pi^0_1$-definable subsets of natural numbers, while coanalytic sets correspond to $\Sigma^0_1$-definable subsets of natural numbers (cf. \cite{addison1954a} p. 80). The explanation is simple: Addison notes that Lusin himself had suggested that analytic sets have many of the same properties as closed sets, while co-analytic sets have many of the same properties as open sets, particularly as concerns their separation properties.\footnote{Indeed, it seems that, for Lusin, the term ``separation'' stems from the $T_1$-axiom in topology, which says that two distinct points can be ``separated'' by disjoint open sets (\cite{lusin1930b} p. 57); elsewhere he continues the topological analogies and suggests that separation principles give us a qualitative notion of distance (\cite{lusin1930a} pp. 65-66).} But on Addison's effectivization of the Borel sets, one has the effectively closed sets $K_1^-$ are precisely the $\Pi^0_1$-definable classes, while the effectively open sets $K_1^+$ are precisely the $\Sigma^0_1$-definable classes.

\subsection{$\mathsf{ATR}_0$ as a formal system}\label{sec:ATRformal}

Nearly simultaneous to Addison's \emph{effectivization} of descriptive set theory, Kreisel was examining which theorems of classical descriptive set theory were \emph{provable} from the systems of predicative analysis, such as $\mathsf{\Sigma^1_1\mbox{-}AC}_0$, which was mentioned above in \S\ref{subsec:Konpreddefn}. He showed that the perfect set theorem for closed sets failed in these systems. In particular, if we restrict attention to closed subsets of Cantor space, these can be represented as paths through infinite binary branching trees. \citet{kreisel1959b} basically showed that there was an effectively closed set whose decomposition into a perfect set and a countable set was such that neither part of the decomposition was hyperarithmetic (cf. \cite{cenzer2012a} \S{IV.7}).\footnote{The qualifier ``basically'' is due only to the fact that \citet{kreisel1959b} works with closed subsets of the unit interval.} 

The traditional proof of the theorem that a non-empty closed set decomposes into a countable set and a perfect set explicitly uses ordinals.\footnote{Cf. \cite{kechris1995a} p. 33, \cite{cenzer2012a} \S{V1}. It turns out that this theorem reverses to $\mathsf{\Pi^1_1\mbox{-}CA}_0$ (cf. \cite{simpson2009} pp. 219-220).} For a closed set $C$, one defines a decreasing sequence of closed sets $C_{\alpha}$ by setting $C_{0}=C$, taking intersections at limits, and setting $C_{\alpha+1}$ to be the non-isolated points of $C_{\alpha}$. Hence it was also natural for Kreisel to study facts about ordinals and well-orderings in these systems. In 1963, Kreisel further showed that these systems did not prove that any two well-orders were comparable, in that they did not prove that one was isomorphic to an initial segment of the other. This result was reported, and a proof given, in \cite{harrison1968a} (cf. pp. 531--532).

The approach that Harrison adopted to these problems was to consider analogues of the hyperarithmetic hierarchy. This is the hierarchy of subsets~$H_a$ of natural numbers obtained by iterating the Turing jump along codes~$a$ for computable ordinals, and it was a result of Kleene's that the $\Delta^1_1$-definable sets were precisely those computable from some element of this hierarchy (cf. \cite{kleene1955a}, \cite{sacks1990a} pp. 24, 31, \cite{ash2000a} p. 81). Following \cite{feferman1962ac},  Harrison considered `pseudo-hierarchies' $H^{\ast}_a$ which were defined just like $H_a$ but with the exception that the codes~$a$ were from a linear order which, while not truly a well-order, at least had no hyperarithmetic infinite descending sequences. Harrison used this to show that one could not prove, in the systems Kreisel had considered, that for every code~$a$ for a computable ordinal that the set $H_a$ existed (cf. \cite{harrison1968a} pp. 536, 542).

Friedman's dissertation interacts with these results in several different ways. First, he showed that one of Kreisel's systems ($\mathsf{\Sigma^1_1\mbox{-}DC}$) is $\Pi^1_2$-conservative over another ($\mathsf{\Delta^1_1\mbox{-}CA}$), by showing that the existence of these hierarchies could go proxy for the comparability of well-orderings.\footnote{Cf. \cite{friedman1967a} p. 10. It should be noted that the theories which Friedman's results concern include full induction. There is no conservation result for the systems without full induction. Cf. \cite{simpson2009} p. 347, p. 381.} Second, he showed that it was not the case that for every pseudo-code $a$ that the set $H^{\ast}_a$ exists (cf. \cite{friedman1967a} pp. 13-14, cf. \cite{harrison1968a} p. 542).

Against this background, it can retrospectively appear almost necessary to consider axiomatic renditions of such hierarchies. Since these hierarchies are defined by transfinite recursion on ordinals, this thus recommends an axiom to the effect that one can effect transfinite recursion along well-orders. And indeed, the principle $\mathsf{ATR_0}$ expresses that if $\theta$ is any arithmetical operator and if $\alpha$ is a code for a well-order, then the effective union $\oplus_{\beta<\alpha} \theta_{\beta}$ exists, where this is defined recursively by
\begin{equation*}
\theta_0 = \emptyset, \hspace{10mm} \theta_{\beta+1} = \theta(\theta_{\beta}), \hspace{10mm} \oplus_{\beta<\gamma} \theta_{\beta}, \mbox{ if $\gamma<\alpha$ limit}
\end{equation*}
This axiom was first articulated by Friedman in his 1974 address, where he mentions that it is equivalent to the comparability of countable well-orders (cf. \cite{friedman1975} p. 240, \cite{friedman1976} p. 559). In Simpson's monograph, it is further shown that $\mathsf{ATR_0}$ is equivalent over $\mathsf{ACA}_0$ to the ``oracle-version'' of the claim that $H_a$ exists for any code $a$ for a computable well-order (\cite{simpson2009} Theorem VIII.3.15 p. 328), as well as ``every uncountable analytic set has a perfect subset'' (cf. \cite{simpson2009} Theorem V.5.5 p. 193).

This axiomatization of $\mathsf{ATR_0}$ accords well with Borel and Lusin's idea that an ordinal was an ``abbreviated notation'' for ``the order in which must be effected a countable number of operations'': for, $\mathsf{ATR_0}$ directly postulates an ability to do recursion along such an ordinal. Further, since this is how ordinals are actually used in classical descriptive set theory-- namely as indices for certain stages in a construction -- it is not surprising $\mathsf{ATR_0}$ ends up reversing to many statements of descriptive set theory. The first of these reversals were from Steel's dissertation, written under Addison and Simpson. In particular, Steel showed that $\mathsf{ATR_0}$ was equivalent to weak forms of determinacy, namely to the case where the winning set was an effective open set (cf. \cite{steel1977} p. 15, \cite{simpson2009} p. 208). While this result is now a classic of reverse mathematics, its constituent notions -- the axiom $\mathsf{ATR}_0$ itself and the effectivization of the Borel hierarchy -- were of a long genesis, whose steps we have sought to record and document.\footnote{Of course, other results from effective descriptive set theory, such as Silver's Theorem and slightly stronger forms of determinacy, reverse to systems stronger than $\mathsf{ATR}_0$ (cf. \cite{simpson2009} pp. 229, 235). Further, yet stronger forms of determinacy are not even provable in full second-order Peano arithmetic (cf. \cite{friedman1970}, \cite{montalban2012}).}

Finally, it is worth mentioning one further aspect of the history of the axiom $\mathsf{ATR_0}$. Again, this axiom says that the we can do recursion along a well-order. But there are different ways of making the result of a computation available. For instance, instead of making the entire computation process available, one might simply make a certain outcome of the result available. In the case where the governing operator is monotonic, it is of course well-known that there will exist a least-fixed point, and so one might rather postulate axioms asserting the existence of this point. This was the impetus behind the many theories of inductive definitions studied by Kreisel, Feferman and others, and surveyed in \cite{feferman1981a}, \cite{feferman1981b}. One of the motivations for this study came from descriptive set theory: ``Inductive definitions formulate rules for generating mathematical objects. The process of inductive generation is used frequently in mathematics; for example, it is used to obtain the subgroup of a group $G$ (generated from a given subset of $G$) or the Borel sets of a topological space'' (\cite{feferman1981b} p. 18). 

\section{Conclusions}\label{sec:conclusions}

In the previous pages we have set out in detail the history of the constituent subsystems $\mathsf{RCA}_0$, $\mathsf{WKL}_0$, $\mathsf{ACA}_0$, $\mathsf{ATR}_0$ of second-order arithmetic. We want to close with a few reflections, of a more general character, on what this suggests for our understanding of reverse mathematics and the many parts of the history of logic with which it interacts.

Let us begin with reverse mathematics. Our history suggests that there is a long-standing tradition, within the study of the subsystems, of evaluating candidate axioms by studying the implications between these axioms and principles of a more ostenstibly mathematical character. For instance, we saw this early in the debate between Poincar\'e and Zermelo in \S\ref{subsec:RPZWsetexistence}, where they both accepted that the predicative perspective must prove the Fundamental Theorem of Algebra if this perspective was to be acceptable. Likewise, in \S\ref{sec:basisthm} we saw the how the Completeness Theorem's requirement of the existence of non-computable sets casts doubt as to its constructive credentials. Episodes such as these then suggest that the organization of the study of subsystems of second-order arithmetic around reversals simply explicitly thematizes this long-standing element of the practice.

Further, the history of this subject suggests an alternative to the received view on the significance of reversals. The received view, due to Simpson and repeated in nearly every talk and paper on reverse mathematics, is that reversals are significant because they measure the set-existence principles implicit in ordinary mathematics.\footnote{See  \citet{simpson1988a}, \citet[I.1, I.9]{simpson2009}. But compare \citet[\S{3}]{friedman2000}, where it is suggested that one might considering replacing the `mutual derivability' that is characteristic of reversals with the related notion of `mutual interpretability'; likewise, Drake had at one point suggested that reversals might be interesting because they track `consistency strength' (\citet[\S{2.4}]{drake1989a}). See \citet{walsh2014a} for a discussion of what epistemic notions may or may not be tracked by interpretability.} But if ``set existence principle'' means ``instance of the comprehension schema,'' then it leaves out $\mathsf{WKL}_0$ and $\mathsf{ATR}_0$. If in response to this, one broadens the definition to include any sentence beginning with ``$\forall \; X \; \exists \; Y \ldots$'', then by trivial syntactic manipulations every sentence can be made to be equivalent to a set-existence principle, and then measuring set-existence would be just the same as sorting out the very fine-grained equivalence classes of mutual derivability.

A deflationary alternative, suggested by our history, is that reversals are significant simply because the axioms of the subsystems are antecedently significant, and showing something to be equivalent to such a subsystem provides it with additional meaning and significance. In the sections devoted to $\mathsf{ACA}_0$, $\mathsf{RCA}_0$, and $\mathsf{ATR}_0$, we described various philosophical viewpoints which can be seen to motivate these positions. The knowledge that a principle  of mathematics is equivalent to one of these subsystems then further broadens our knowledge of this viewpoint: we then know that the principle is justifiable from this perspective. And if too few such important principles are so justified, then we have a reason to move beyond this perspective.

But such a broadening can be effected by things other than a reversal \emph{per se}. In addition to knowing what kinds of algebra and analysis are justifiable from a given perspective, we might want to know what sorts of computation-like processes are justifiable from a given perspective. For instance, knowing whether the perspective is compatible with there only being computable sets, or knowing whether the perspective is compatible with there only being primitive recursive provably total functions, seems to be on a par with knowing what kinds of algebra and analysis may be done within a system. Bishop's constructivism prided itself on being compatible with only computable sets, and Tait's understanding of finitism restricted the available functions to the primitive recursive functions. To understand how far one is from such constructivism and finitism is to understand one's computational resources.\footnote{This is related to the suggestion in \citet[\S{2}]{shore2010} that what is really interesting in reverse mathematics is showing that two principles have the same $\omega$-models. But two theories can have the same $\omega$-models without having the same provably recursive functions, and vice-versa. Indeed, it is this very phenomena which suggested the shift-- mentioned in the previous footnote-- from `mutual derivability' to `mutual interpretability'.}

On this picture, it is then natural to think that there will be certain subsystems which, while not themselves corresponding to any philosophical or foundational viewpoint, serve as guideposts for the calibration of such resources. Our history suggests that $\mathsf{WKL}_0$ is such a subsystem: indeed, it arose as one of the first markers separating mathematics which is compatible with only computable sets from mathematics which requires some non-computable sets. In this sense then, of course a reversal to $\mathsf{WKL}_0$ does tell us something about set existence, namely, it mandates a wide array of non-computable sets. But saying this is not to say that there is a univocal notion of set existence which captures both $\mathsf{WKL}_0$ and the subsystems formed from restrictions on comprehension schema. Rather, the thought would simply be: to understand a subsystem and the foundational viewpoint which it represents is to understand the mathematics which it is consistent with, and computation is just as much a part of mathematics as algebra and analysis.

But there is much in our history that ought be of interest to those without prior interests in reverse mathematics as such. Hence we want to close by highlighting a distinctive feature of this history, namely: the crucial intermediary role that certain figures play in transmitting results between different areas of mathematical logic and different foundational enterprises. While not having the stature of a Hilbert or G\"odel or Turing, figures such as Kleene and Kreisel played a distinctive role as instigators of interdisciplinarity within logic, a role which is omitted from the usual description of their achievements. For instance, while Kleene is rightfully regarded as one of the founders of the theory of computation, in our brief history he also plays an important role in forging interactions between computability theory and other sub-fields of logic. To briefly recapitulate: (1) in \S\ref{subsec:GMKeffectiveanalysis} we mentioned the startling fact that Kleene's characterization of the computable sets as the $\Delta^0_1$-definable sets was inspired by Souslin's theorem, (2) in \S\ref{sec:basisthm} we described how the study of the effective properties of K\"onig's Lemma was initiated by Kleene as a study of Brouwer's constructive foundations of analysis, (3) and in \S\ref{subsec:borel} we described how it was Addison's dissertation, written under Kleene, that transformed Borel's attempt to study ``calculable'' restrictions within analysis into present-day effective descriptive set theory. 

Likewise, while today Kreisel is primarily known among philosophers for his claims about second-order logic and the determinacy of the continuum hypothesis, and among mathematicians for his basis theorems (cf. \cite{kreisel1967}, \cite{cenzer2012a} \S{III.2}), our history shows that he played an important role in translating and collating results between different areas and different communities. Again, to summarize: (1) in \S\S\ref{subsec:Konpreddefn}-\ref{sec:KWFpredprov} we mentioned Kreisel's role in emphasizing the important difference between predicative provability and definability, (2) in \S\ref{sec:basisthm} we described his role in bringing computability-theoretic results to bear on the status of the completeness theorem within the constructivist framework, (3) and in \S\ref{sec:ATRformal} we described how Kreisel's work on the computational strength of the perfect set theorem for closed sets led to Harrison's investigations, which were the immediate antecedent of some of Friedman's results in his dissertation.

If one were writing merely the history of one branch of mathematical logic or foundational enterprise in isolation from the rest, it would be difficult to appreciate the combined magnitude of these contributions of Kleene and Kreisel to the history of logic in the last century. No doubt we have not said all there is to say about their contributions or that of other figures such as Feferman or Friedman whose cross-disciplinary work is better known.  None of these figures adhered to or is associated with a specific foundational standpoint.  But our history has also served to illustrate how their extension of results and methods which grew out of these standpoints served to guide the development of mathematical logic in the second half of the last century -- albeit in indirect and occasionally surprising directions.

\section{Acknowledgments}

This paper has been measurably improved by us having had the opportunity to present versions of it at the following events: Computability Theory and Foundations of Mathematics at the Tokyo Institute of Technology on September 8, 2015, the Logic Seminar at the University of California, Irvine on January 13, 2016; at a special session on the history and philosophy of logic at the North American annual meeting of the Association for Symbolic Logic on May 24, 2016; and at the  Southern California History and Philosophy of Logic and Mathematics Group on June 15, 2016. Thanks is owed in particular to the following people who gave us invaluable feedback and comments on these and other occasions:  Andrew Arana, Kyle Banick, Patrica Blanchette, Tim Button, Michael Detlefsen, Benedict Eastaugh, Samuel Eklund, Ulrich Kohlenbach, Greg Lauro, Adam Harmer, Jeremy Heis, Robert Lubarsky, Richard Mendelsohn,  Christopher Mitsch, Erich Reck, Gillian Russell, Stewart Shapiro, Stephen Simpson, Jeff Schatz, Will Stafford, Clinton Tolley, Alasdair Urquhart, Kai Wehmeier, Wilfried Sieg, Kino Zhao, and Richard Zach. 

\bibliography{prehistory.bib}

\end{document}